\documentclass{amsart}
\usepackage{fullpage}

\def \R{\mathbb{R}}
\def \Z{\mathbb{Z}}

\usepackage{color}
\usepackage{amsfonts}
\usepackage{amsmath, amscd, amssymb, tikz-cd}

\newcommand{\g}{\mathfrak{g}}
\newcommand{\h}{\mathfrak{h}}
\newcommand{\el}{\mathfrak{l}}
\newcommand{\Q}{\mathbb{Q}}
\newcommand{\Prim}{\mathrm{Prim}}
\newcommand{\U}{\mathrm{U}}
\newcommand{\m}{\mathfrak{m}}
\newcommand{\Nalt}{\mathrm{N_{\mathrm{alt}}}}
\newcommand{\Aut}{\mathrm{Aut}}
\newcommand{\LNalt}{\mathrm{LN_{\mathrm{alt}}}}
\newcommand{\Spin}{\mathrm{Spin}}

\input{xy}
\xyoption{all}

\newcount\grcalca
 \newcount\grcalcb
 \newcount\grcalcc
 \newcount\grcalcd
 \newcount\grcalce
 \newcount\grcalcf
 \newcount\grcalcg
 \newcount\grcalch
 \newcount\grcalci
 \newcount\grrow
 \newcount\grcolumn
 \newcount\grwidth
 \newcount\factor
 \newcount\hfactor
 \newcount\qfactor
 \newcount\tfactor
 \newcount\sfactor
 \newcount\dfactor
 \hfactor = \factor
 \divide    \hfactor by 2
 \qfactor = \hfactor
 \divide    \qfactor by 2
 \tfactor = \factor
 \divide    \tfactor by 3
 \sfactor = \factor
 \divide    \sfactor by 6
 \dfactor = \factor
 \divide    \dfactor by 12

 \newcommand{\gbeg}[2]{
   \unitlength=1pt
   \grrow = #2
   \grcolumn = 0
   \grcalca = #1
   \grcalcb = #2
   \multiply \grcalca by \factor
   \grwidth = \grcalca
   \multiply \grcalcb by \factor
   \begin{minipage}{\grcalca pt}
   \begin{picture}(\grcalca,\grcalcb)
   \advance \grcalcb by -\factor
   }

 \newcommand{\gend}{
   \end{picture}
   {\vskip2.5ex}
   \end{minipage} }

 \newcommand{\gnl}{
   \advance \grrow by -1
   \grcolumn = 0}

 \newcommand{\gvac}[1]{
   \advance \grcolumn by #1}

 \newcommand{\gcl}[1]{
   \grcalca = \grcolumn
   \multiply \grcalca by \factor
   \advance \grcalca by \hfactor
   \grcalcb = \grrow
   \multiply \grcalcb by \factor
   \grcalcc = #1
   \multiply \grcalcc by \factor
   \put(\grcalca,\grcalcb) {\line(0,-1){\grcalcc}}
   \advance \grcolumn by 1}

 \newcommand{\gcn}[4]{
   \grcalca = \grcolumn
   \multiply \grcalca by \factor
   \grcalci = #3
   \multiply \grcalci by \hfactor
   \advance \grcalca by \grcalci
   \grcalcb = \grcolumn
   \multiply \grcalcb by \factor
   \grcalci = #3
   \advance \grcalci by #4
   \multiply \grcalci by \qfactor
   \advance \grcalcb by \grcalci
   \grcalcc = \grcolumn
   \multiply \grcalcc by \factor
   \grcalci = #4
   \multiply \grcalci by \hfactor
   \advance \grcalcc by \grcalci
   \grcalcd = \grrow
   \multiply \grcalcd by \factor
   \grcalce = \grrow
   \multiply \grcalce by \factor
   \grcalci = #2
   \multiply \grcalci by \tfactor
   \advance \grcalce by -\grcalci
   \grcalcf = \grrow
   \multiply \grcalcf by \factor
   \grcalci = #2
   \multiply \grcalci by \hfactor
   \advance \grcalcf by -\grcalci
   \grcalcg = \grrow
   \multiply \grcalcg by \factor
   \grcalci = #2
   \multiply \grcalci by \tfactor
   \multiply \grcalci by 2
   \advance \grcalcg by -\grcalci
   \grcalch = \grrow
   \advance \grcalch by -#2
   \multiply \grcalch by \factor
   \qbezier(\grcalca,\grcalcd)(\grcalca,\grcalce)(\grcalcb,\grcalcf)
   \qbezier(\grcalcb,\grcalcf)(\grcalcc,\grcalcg)(\grcalcc,\grcalch)
   \advance \grcolumn by #1}

 \newcommand{\gnot}[1]{
   \grcalca = \grcolumn
   \multiply \grcalca by \factor
   \advance \grcalca by \hfactor
   \grcalcb = \grrow
   \multiply \grcalcb by \factor
   \advance \grcalcb by -\hfactor
   \put(\grcalca,\grcalcb) {\makebox(0,0){$\scriptstyle #1$}} }

 \newcommand{\got}[2]{
   \grcalca = \grcolumn
   \multiply \grcalca by \factor
   \grcalcc = #1
   \multiply \grcalcc by \hfactor
   \advance \grcalca by \grcalcc
   \grcalcb = \grrow
   \multiply \grcalcb by \factor
   \advance \grcalcb by -\tfactor
   \advance \grcalcb by -\tfactor
   \put(\grcalca,\grcalcb){\makebox(0,0)[b]{$#2$}}
   \advance \grcolumn by #1}

 \newcommand{\gob}[2]{
   \grcalca = \grcolumn
   \multiply \grcalca by \factor
   \grcalcc = #1
   \multiply \grcalcc by \hfactor
   \advance \grcalca by \grcalcc
   \put(\grcalca,0){\makebox(0,0)[b]{$#2$}}
   \advance \grcolumn by #1}

 \newcommand{\gmu}{
   \grcalca = \grcolumn
   \advance \grcalca by 1
   \multiply \grcalca by \factor
   \grcalcb = \grrow
   \multiply \grcalcb by \factor
   \grcalcc = \factor
   \advance \grcalcc by \hfactor
   \put(\grcalca,\grcalcb){\oval(\factor,\grcalcc)[b]}
   \advance \grcalcb by -\hfactor
   \advance \grcalcb by -\qfactor
   \put(\grcalca,\grcalcb) {\line(0,-1){\qfactor}}
   \advance \grcolumn by 2}

 \newcommand{\gcmu}{
   \grcalca = \grcolumn
   \advance \grcalca by 1
   \multiply \grcalca by \factor
   \grcalcb = \grrow
   \advance \grcalcb by -1
   \multiply \grcalcb by \factor
   \grcalcc = \factor
   \advance \grcalcc by \hfactor
   \put(\grcalca,\grcalcb){\oval(\factor,\grcalcc)[t]}
   \advance \grcalcb by \factor
   \put(\grcalca,\grcalcb) {\line(0,-1){\qfactor}}
   \advance \grcolumn by 2}

 \newcommand{\glm}{
   \grcalca = \grcolumn
   \multiply \grcalca by \factor
   \advance \grcalca by \hfactor
   \grcalcb = \grcalca
   \advance \grcalcb by \factor
   \grcalcc = \grrow
   \multiply \grcalcc by \factor
   \grcalcd = \grcalcc
   \advance \grcalcd by -\tfactor
   \grcalce = \grcalcd
   \advance \grcalce by -\tfactor
   \put(\grcalca, \grcalcc){\line(0,-1){\tfactor}}
   \put(\grcalca, \grcalcd){\line(1,0){\factor}}
   \put(\grcalca, \grcalcd){\line(3,-1){\factor}}
   \put(\grcalcb, \grcalcc){\line(0,-1){\factor}}
   \advance \grcolumn by 2}

 \newcommand{\grm}{
   \grcalcb = \grcolumn
   \multiply \grcalcb by \factor
   \advance \grcalcb by \hfactor
   \grcalca = \grcalcb
   \advance \grcalca by \factor
   \grcalcc = \grrow
   \multiply \grcalcc by \factor
   \grcalcd = \grcalcc
   \advance \grcalcd by -\tfactor
   \grcalce = \grcalcd
   \advance \grcalce by -\tfactor
   \put(\grcalca, \grcalcc){\line(0,-1){\tfactor}}
   \put(\grcalca, \grcalcd){\line(-1,0){\factor}}
   \put(\grcalca, \grcalcd){\line(-3,-1){\factor}}
   \put(\grcalcb, \grcalcc){\line(0,-1){\factor}}
   \advance \grcolumn by 2}

 \newcommand{\glcm}{
   \grcalca = \grcolumn
   \multiply \grcalca by \factor
   \advance \grcalca by \hfactor
   \grcalcb = \grcalca
   \advance \grcalcb by \factor
   \grcalcc = \grrow
   \advance \grcalcc by -1
   \multiply \grcalcc by \factor
   \grcalcd = \grcalcc
   \advance \grcalcd by \tfactor
   \grcalce = \grcalcd
   \advance \grcalce by \tfactor
   \put(\grcalca, \grcalcc){\line(0,1){\tfactor}}
   \put(\grcalca, \grcalcd){\line(1,0){\factor}}
   \put(\grcalca, \grcalcd){\line(3,1){\factor}}
   \put(\grcalcb, \grcalcc){\line(0,1){\factor}}
   \advance \grcolumn by 2}

 \newcommand{\grcm}{
   \grcalcb = \grcolumn
   \multiply \grcalcb by \factor
   \advance \grcalcb by \hfactor
   \grcalca = \grcalcb
   \advance \grcalca by \factor
   \grcalcc = \grrow
   \advance \grcalcc by -1
   \multiply \grcalcc by \factor
   \grcalcd = \grcalcc
   \advance \grcalcd by \tfactor
   \grcalce = \grcalcd
   \advance \grcalce by \tfactor
   \put(\grcalca, \grcalcc){\line(0,1){\tfactor}}
   \put(\grcalca, \grcalcd){\line(-1,0){\factor}}
   \put(\grcalca, \grcalcd){\line(-3,1){\factor}}
   \put(\grcalcb, \grcalcc){\line(0,1){\factor}}
   \advance \grcolumn by 2}

 \newcommand{\gwmu}[1]{
   \grcalca = \grcolumn
   \multiply \grcalca by \factor
   \grcalcd = \hfactor
   \multiply \grcalcd by #1
   \advance \grcalca by \grcalcd
   \grcalcb = \grrow
   \multiply \grcalcb by \factor
   \grcalcc = \factor
   \advance \grcalcc by \hfactor
   \grcalcd = #1
   \advance \grcalcd by -1
   \multiply \grcalcd by \factor
   \put(\grcalca,\grcalcb){\oval(\grcalcd,\grcalcc)[b]}
   \advance \grcalcb by -\hfactor
   \advance \grcalcb by -\qfactor
   \put(\grcalca,\grcalcb) {\line(0,-1){\qfactor}}
   \advance \grcolumn by #1}

 \newcommand{\gwcm}[1]{
   \grcalca = \grcolumn
   \multiply \grcalca by \factor
   \grcalcd = \hfactor
   \multiply \grcalcd by #1
   \advance \grcalca by \grcalcd
   \grcalcb = \grrow
   \advance \grcalcb by -1
   \multiply \grcalcb by \factor
   \grcalcc = \factor
   \advance \grcalcc by \hfactor
   \grcalcd = #1
   \advance \grcalcd by -1
   \multiply \grcalcd by \factor
   \put(\grcalca,\grcalcb){\oval(\grcalcd,\grcalcc)[t]}
   \advance \grcalcb by \factor
   \put(\grcalca,\grcalcb) {\line(0,-1){\qfactor}}
   \advance \grcolumn by #1}

 \newcommand{\gwmuc}[1]{
   \grcalca = \grcolumn
   \multiply \grcalca by \factor
   \advance \grcalca by \hfactor
   \grcalcb = \grrow
   \multiply \grcalcb by \factor
   \grcalcc = #1
   \advance \grcalcc by -1
   \multiply \grcalcc by \factor
   \put(\grcalca,\grcalcb){\line(1,0){\grcalcc}}
   \advance \grcalca by -\hfactor
   \grcalcd = \hfactor
   \multiply \grcalcd by #1
   \advance \grcalca by \grcalcd
   \grcalcc = \factor
   \advance \grcalcc by \hfactor
   \grcalcd = #1
   \advance \grcalcd by -1
   \multiply \grcalcd by \factor
   \put(\grcalca,\grcalcb){\oval(\grcalcd,\grcalcc)[b]}
   \advance \grcalcb by -\hfactor
   \advance \grcalcb by -\qfactor
   \put(\grcalca,\grcalcb) {\line(0,-1){\qfactor}}
   \advance \grcolumn by #1}

 \newcommand{\gwcmc}[1]{
   \grcalca = \grcolumn
   \multiply \grcalca by \factor
   \advance \grcalca by \hfactor
   \grcalcb = \grrow
   \multiply \grcalcb by \factor
   \advance \grcalcb by -\factor
   \grcalcc = #1
   \advance \grcalcc by -1
   \multiply \grcalcc by \factor
   \put(\grcalca,\grcalcb){\line(1,0){\grcalcc}}
   \grcalcd = #1
   \advance \grcalcd by -1
   \multiply \grcalcd by \hfactor
   \advance \grcalca by \grcalcd
   \grcalcc = \factor
   \advance \grcalcc by \hfactor
   \grcalcd = #1
   \advance \grcalcd by -1
   \multiply \grcalcd by \factor
   \put(\grcalca,\grcalcb){\oval(\grcalcd,\grcalcc)[t]}
   \advance \grcalcb by \factor
   \put(\grcalca,\grcalcb) {\line(0,-1){\qfactor}}
   \advance \grcolumn by #1}

 \newcommand{\gev}{
   \grcalca = \grcolumn
   \advance \grcalca by 1
   \multiply \grcalca by \factor
   \grcalcb = \grrow
   \multiply \grcalcb by \factor
   \grcalcc = \factor
   \advance \grcalcc by \hfactor
   \put(\grcalca,\grcalcb){\oval(\factor,\grcalcc)[b]}
   \advance \grcolumn by 2}

 \newcommand{\gdb}{
   \grcalca = \grcolumn
   \advance \grcalca by 1
   \multiply \grcalca by \factor
   \grcalcb = \grrow
   \advance \grcalcb by -1
   \multiply \grcalcb by \factor
   \grcalcc = \factor
   \advance \grcalcc by \hfactor
   \put(\grcalca,\grcalcb){\oval(\factor,\grcalcc)[t]}
   \advance \grcolumn by 2}

 \newcommand{\gwev}[1]{
   \grcalca = \grcolumn
   \multiply \grcalca by \factor
   \grcalcd = \hfactor
   \multiply \grcalcd by #1
   \advance \grcalca by \grcalcd
   \grcalcb = \grrow
   \multiply \grcalcb by \factor
   \grcalcc = \factor
   \advance \grcalcc by \hfactor
   \grcalcd = #1
   \advance \grcalcd by -1
   \multiply \grcalcd by \factor
   \put(\grcalca,\grcalcb){\oval(\grcalcd,\grcalcc)[b]}
   \advance \grcolumn by #1}

 \newcommand{\gwdb}[1]{
   \grcalca = \grcolumn
   \multiply \grcalca by \factor
   \grcalcd = \hfactor
   \multiply \grcalcd by #1
   \advance \grcalca by \grcalcd
   \grcalcb = \grrow
   \advance \grcalcb by -1
   \multiply \grcalcb by \factor
   \grcalcc = \factor
   \advance \grcalcc by \hfactor
   \grcalcd = #1
   \advance \grcalcd by -1
   \multiply \grcalcd by \factor
   \put(\grcalca,\grcalcb){\oval(\grcalcd,\grcalcc)[t]}
   \advance \grcolumn by #1}

 \newcommand{\gbr}{
   \grcalca = \grcolumn
   \multiply \grcalca by \factor
   \advance \grcalca by \hfactor
   \grcalcb = \grcalca
   \advance \grcalcb by \hfactor
   \grcalcc = \grcalca
   \advance \grcalcc by \factor
   \grcalcd = \grrow
   \multiply \grcalcd by \factor
   \grcalce = \grcalcd
   \advance \grcalce by -\tfactor
   \grcalcf = \grcalcd
   \advance \grcalcf by -\hfactor
   \grcalcg = \grcalce
   \advance \grcalcg by -\tfactor
   \grcalch = \grcalcd
   \advance \grcalch by -\factor
   \qbezier(\grcalca,\grcalcd)(\grcalca,\grcalce)(\grcalcb,\grcalcf)
   \qbezier(\grcalcb,\grcalcf)(\grcalcc,\grcalcg)(\grcalcc,\grcalch)
   \advance \grcalcf by -\dfactor
   \advance \grcalcb by -\sfactor
   \qbezier(\grcalca,\grcalch)(\grcalca,\grcalcg)(\grcalcb,\grcalcf)
   \advance \grcalcf by \sfactor
   \advance \grcalcb by \tfactor
   \qbezier(\grcalcc,\grcalcd)(\grcalcc,\grcalce)(\grcalcb,\grcalcf)
   \advance \grcolumn by 2}

 \newcommand{\gibr}{
   \grcalca = \grcolumn
   \multiply \grcalca by \factor
   \advance \grcalca by \hfactor
   \grcalcb = \grcalca
   \advance \grcalcb by \hfactor
   \grcalcc = \grcalca
   \advance \grcalcc by \factor
   \grcalcd = \grrow
   \multiply \grcalcd by \factor
   \grcalce = \grcalcd
   \advance \grcalce by -\tfactor
   \grcalcf = \grcalcd
   \advance \grcalcf by -\hfactor
   \grcalcg = \grcalce
   \advance \grcalcg by -\tfactor
   \grcalch = \grcalcd
   \advance \grcalch by -\factor
   \qbezier(\grcalcc,\grcalcd)(\grcalcc,\grcalce)(\grcalcb,\grcalcf)
   \qbezier(\grcalcb,\grcalcf)(\grcalca,\grcalcg)(\grcalca,\grcalch)
   \advance \grcalcf by -\dfactor
   \advance \grcalcb by \sfactor
   \qbezier(\grcalcc,\grcalch)(\grcalcc,\grcalcg)(\grcalcb,\grcalcf)
   \advance \grcalcf by \sfactor
   \advance \grcalcb by -\tfactor
   \qbezier(\grcalca,\grcalcd)(\grcalca,\grcalce)(\grcalcb,\grcalcf)
   \advance \grcolumn by 2}

 \newcommand{\gbrc}{
   \grcalca = \grcolumn
   \multiply \grcalca by \factor
   \advance \grcalca by \hfactor
   \grcalcb = \grcalca
   \advance \grcalcb by \hfactor
   \grcalcc = \grcalca
   \advance \grcalcc by \factor
   \grcalcd = \grrow
   \multiply \grcalcd by \factor
   \grcalce = \grcalcd
   \advance \grcalce by -\tfactor
   \grcalcf = \grcalcd
   \advance \grcalcf by -\hfactor
   \grcalcg = \grcalce
   \advance \grcalcg by -\tfactor
   \grcalch = \grcalcd
   \advance \grcalch by -\factor
   \put(\grcalcb,\grcalcf){\circle{\hfactor}}
   \qbezier(\grcalca,\grcalcd)(\grcalca,\grcalce)(\grcalcb,\grcalcf)
   \qbezier(\grcalcb,\grcalcf)(\grcalcc,\grcalcg)(\grcalcc,\grcalch)
   \advance \grcalcf by -\dfactor
   \advance \grcalcb by -\sfactor
   \qbezier(\grcalca,\grcalch)(\grcalca,\grcalcg)(\grcalcb,\grcalcf)
   \advance \grcalcf by \sfactor
   \advance \grcalcb by \tfactor
   \qbezier(\grcalcc,\grcalcd)(\grcalcc,\grcalce)(\grcalcb,\grcalcf)
   \advance \grcolumn by 2}

 \newcommand{\gibrc}{
   \grcalca = \grcolumn
   \multiply \grcalca by \factor
   \advance \grcalca by \hfactor
   \grcalcb = \grcalca
   \advance \grcalcb by \hfactor
   \grcalcc = \grcalca
   \advance \grcalcc by \factor
   \grcalcd = \grrow
   \multiply \grcalcd by \factor
   \grcalce = \grcalcd
   \advance \grcalce by -\tfactor
   \grcalcf = \grcalcd
   \advance \grcalcf by -\hfactor
   \grcalcg = \grcalce
   \advance \grcalcg by -\tfactor
   \grcalch = \grcalcd
   \advance \grcalch by -\factor
   \put(\grcalcb,\grcalcf){\circle{\hfactor}}
   \qbezier(\grcalcc,\grcalcd)(\grcalcc,\grcalce)(\grcalcb,\grcalcf)
   \qbezier(\grcalcb,\grcalcf)(\grcalca,\grcalcg)(\grcalca,\grcalch)
   \advance \grcalcf by -\dfactor
   \advance \grcalcb by \sfactor
   \qbezier(\grcalcc,\grcalch)(\grcalcc,\grcalcg)(\grcalcb,\grcalcf)
   \advance \grcalcf by \sfactor
   \advance \grcalcb by -\tfactor
   \qbezier(\grcalca,\grcalcd)(\grcalca,\grcalce)(\grcalcb,\grcalcf)
   \advance \grcolumn by 2}

 \newcommand{\gu}[1]{
   \grcalca = \grcolumn
   \multiply \grcalca by \factor
   \grcalcd = \hfactor
   \multiply \grcalcd by #1
   \advance \grcalca by \grcalcd
   \grcalcb = \grrow
   \advance \grcalcb by -1
   \multiply \grcalcb by \factor
   \put(\grcalca,\grcalcb) {\line(0,1){\hfactor}}
   \advance \grcalcb by \hfactor
   \put(\grcalca,\grcalcb) {\circle*{3}}
   \advance \grcolumn by #1}

 \newcommand{\gcu}[1]{
   \grcalca = \grcolumn
   \multiply \grcalca by \factor
   \grcalcd = \hfactor
   \multiply \grcalcd by #1
   \advance \grcalca by \grcalcd
   \grcalcb = \grrow
   \multiply \grcalcb by \factor
   \put(\grcalca,\grcalcb) {\line(0,-1){\hfactor}}
   \advance \grcalcb by -\hfactor
   \put(\grcalca,\grcalcb) {\circle*{3}}
   \advance \grcolumn by #1}

 \newcommand{\gmp}[1]{
   \grcalca = \grcolumn
   \multiply \grcalca by \factor
   \advance \grcalca by \hfactor
   \grcalcb = \grrow
   \multiply \grcalcb by \factor
   \put(\grcalca,\grcalcb) {\line(0,-1){\dfactor}}
   \advance \grcalcb by -\factor
   \put(\grcalca,\grcalcb) {\line(0,1){\dfactor}}
   \advance \grcalcb by \hfactor
   \grcalcc = \factor
   \advance \grcalcc by -\qfactor
   \put(\grcalca,\grcalcb) {\circle{\grcalcc}}
   \put(\grcalca,\grcalcb) {\makebox(0,0){$\scriptstyle #1$}}
   \advance \grcolumn by 1}

 \newcommand{\gbmp}[1]{
   \grcalca = \grcolumn
   \multiply \grcalca by \factor
   \advance \grcalca by \hfactor
   \grcalcb = \grrow
   \multiply \grcalcb by \factor
   \put(\grcalca,\grcalcb) {\line(0,-1){\dfactor}}
   \advance \grcalcb by -\factor
   \put(\grcalca,\grcalcb) {\line(0,1){\dfactor}}
   \advance \grcalca by -\hfactor
   \advance \grcalca by \dfactor
   \advance \grcalcb by \dfactor
   \grcalcc = \factor
   \advance \grcalcc by -\sfactor
   \put(\grcalca,\grcalcb) {\framebox(\grcalcc,\grcalcc){$\scriptstyle #1$}}
   \advance \grcolumn by 1}

 \newcommand{\gbmpt}[1]{
   \grcalca = \grcolumn
   \multiply \grcalca by \factor
   \advance \grcalca by \hfactor
   \grcalcb = \grrow
   \multiply \grcalcb by \factor
   \put(\grcalca,\grcalcb) {\line(0,-1){\dfactor}}
   \advance \grcalcb by -\factor
   \advance \grcalca by -\hfactor
   \advance \grcalca by \dfactor
   \advance \grcalcb by \dfactor
   \grcalcc = \factor
   \advance \grcalcc by -\sfactor
   \put(\grcalca,\grcalcb) {\framebox(\grcalcc,\grcalcc){$\scriptstyle #1$}}
   \advance \grcolumn by 1}

 \newcommand{\gbmpb}[1]{
   \grcalca = \grcolumn
   \multiply \grcalca by \factor
   \advance \grcalca by \hfactor
   \grcalcb = \grrow
   \multiply \grcalcb by \factor
   \advance \grcalcb by -\factor
   \put(\grcalca,\grcalcb) {\line(0,1){\dfactor}}
   \advance \grcalca by -\hfactor
   \advance \grcalca by \dfactor
   \advance \grcalcb by \dfactor
   \grcalcc = \factor
   \advance \grcalcc by -\sfactor
   \put(\grcalca,\grcalcb) {\framebox(\grcalcc,\grcalcc){$\scriptstyle #1$}}
   \advance \grcolumn by 1}

 \newcommand{\gbmpn}[1]{
   \grcalca = \grcolumn
   \multiply \grcalca by \factor
   \advance \grcalca by \hfactor
   \grcalcb = \grrow
   \multiply \grcalcb by \factor
   \advance \grcalcb by -\factor
   \advance \grcalca by -\hfactor
   \advance \grcalca by \dfactor
   \advance \grcalcb by \dfactor
   \grcalcc = \factor
   \advance \grcalcc by -\sfactor
   \put(\grcalca,\grcalcb) {\framebox(\grcalcc,\grcalcc){$\scriptstyle #1$}}
   \advance \grcolumn by 1}

 \newcommand{\glmptb}{
   \grcalca = \grcolumn
   \multiply \grcalca by \factor
   \advance \grcalca by \hfactor
   \grcalcb = \grrow
   \multiply \grcalcb by \factor
   \put(\grcalca,\grcalcb) {\line(0,-1){\dfactor}}
   \advance \grcalcb by -\factor
   \put(\grcalca,\grcalcb) {\line(0,1){\dfactor}}
   \advance \grcalca by -\hfactor
   \advance \grcalca by \dfactor
   \advance \grcalcb by \dfactor
   \put(\grcalca,\grcalcb) {\line(1,0){\factor}}
   \advance \grcalcb by \factor
   \advance \grcalcb by -\sfactor
   \put(\grcalca,\grcalcb) {\line(1,0){\factor}}
   \grcalcc = \factor
   \advance \grcalcc by -\sfactor
   \put(\grcalca,\grcalcb) {\line(0,-1){\grcalcc}}
   \advance \grcolumn by 1}

 \newcommand{\glmpt}{
   \grcalca = \grcolumn
   \multiply \grcalca by \factor
   \advance \grcalca by \hfactor
   \grcalcb = \grrow
   \multiply \grcalcb by \factor
   \put(\grcalca,\grcalcb) {\line(0,-1){\dfactor}}
   \advance \grcalca by -\hfactor
   \advance \grcalca by \dfactor
   \advance \grcalcb by -\dfactor
   \put(\grcalca,\grcalcb) {\line(1,0){\factor}}
   \advance \grcalcb by -\factor
   \advance \grcalcb by \sfactor
   \put(\grcalca,\grcalcb) {\line(1,0){\factor}}
   \grcalcc = \factor
   \advance \grcalcc by -\sfactor
   \put(\grcalca,\grcalcb) {\line(0,1){\grcalcc}}
   \advance \grcolumn by 1}

 \newcommand{\glmpb}{
   \grcalca = \grcolumn
   \multiply \grcalca by \factor
   \advance \grcalca by \hfactor
   \grcalcb = \grrow
   \multiply \grcalcb by \factor
   \advance \grcalcb by -\factor
   \put(\grcalca,\grcalcb) {\line(0,1){\dfactor}}
   \advance \grcalca by -\hfactor
   \advance \grcalca by \dfactor
   \advance \grcalcb by \dfactor
   \put(\grcalca,\grcalcb) {\line(1,0){\factor}}
   \advance \grcalcb by \factor
   \advance \grcalcb by -\sfactor
   \put(\grcalca,\grcalcb) {\line(1,0){\factor}}
   \grcalcc = \factor
   \advance \grcalcc by -\sfactor
   \put(\grcalca,\grcalcb) {\line(0,-1){\grcalcc}}
   \advance \grcolumn by 1}

 \newcommand{\glmp}{
   \grcalca = \grcolumn
   \multiply \grcalca by \factor
   \advance \grcalca by \dfactor
   \grcalcb = \grrow
   \multiply \grcalcb by \factor
   \advance \grcalcb by -\dfactor
   \put(\grcalca,\grcalcb) {\line(1,0){\factor}}
   \advance \grcalcb by -\factor
   \advance \grcalcb by \sfactor
   \put(\grcalca,\grcalcb) {\line(1,0){\factor}}
   \grcalcc = \factor
   \advance \grcalcc by -\sfactor
   \put(\grcalca,\grcalcb) {\line(0,1){\grcalcc}}
   \advance \grcolumn by 1}

 \newcommand{\gcmptb}{
   \grcalca = \grcolumn
   \multiply \grcalca by \factor
   \advance \grcalca by \hfactor
   \grcalcb = \grrow
   \multiply \grcalcb by \factor
   \put(\grcalca,\grcalcb) {\line(0,-1){\dfactor}}
   \advance \grcalcb by -\factor
   \put(\grcalca,\grcalcb) {\line(0,1){\dfactor}}
   \advance \grcalca by -\hfactor
   \advance \grcalcb by \dfactor
   \put(\grcalca,\grcalcb) {\line(1,0){\factor}}
   \advance \grcalcb by \factor
   \advance \grcalcb by -\sfactor
   \put(\grcalca,\grcalcb) {\line(1,0){\factor}}
   \advance \grcolumn by 1}

 \newcommand{\gcmpt}{
   \grcalca = \grcolumn
   \multiply \grcalca by \factor
   \advance \grcalca by \hfactor
   \grcalcb = \grrow
   \multiply \grcalcb by \factor
   \put(\grcalca,\grcalcb) {\line(0,-1){\dfactor}}
   \advance \grcalcb by -\factor
   \advance \grcalca by -\hfactor
   \advance \grcalcb by \dfactor
   \put(\grcalca,\grcalcb) {\line(1,0){\factor}}
   \advance \grcalcb by \factor
   \advance \grcalcb by -\sfactor
   \put(\grcalca,\grcalcb) {\line(1,0){\factor}}
   \advance \grcolumn by 1}

 \newcommand{\gcmpb}{
   \grcalca = \grcolumn
   \multiply \grcalca by \factor
   \advance \grcalca by \hfactor
   \grcalcb = \grrow
   \multiply \grcalcb by \factor
   \advance \grcalcb by -\factor
   \put(\grcalca,\grcalcb) {\line(0,1){\dfactor}}
   \advance \grcalca by -\hfactor
   \advance \grcalcb by \dfactor
   \put(\grcalca,\grcalcb) {\line(1,0){\factor}}
   \advance \grcalcb by \factor
   \advance \grcalcb by -\sfactor
   \put(\grcalca,\grcalcb) {\line(1,0){\factor}}
   \advance \grcolumn by 1}

 \newcommand{\gcmp}{
   \grcalca = \grcolumn
   \multiply \grcalca by \factor
   \grcalcb = \grrow
   \multiply \grcalcb by \factor
   \advance \grcalcb by -\factor
   \advance \grcalcb by \dfactor
   \put(\grcalca,\grcalcb) {\line(1,0){\factor}}
   \advance \grcalcb by \factor
   \advance \grcalcb by -\sfactor
   \put(\grcalca,\grcalcb) {\line(1,0){\factor}}
   \advance \grcolumn by 1}

 \newcommand{\grmptb}{
   \grcalca = \grcolumn
   \multiply \grcalca by \factor
   \advance \grcalca by \hfactor
   \grcalcb = \grrow
   \multiply \grcalcb by \factor
   \put(\grcalca,\grcalcb) {\line(0,-1){\dfactor}}
   \advance \grcalcb by -\factor
   \put(\grcalca,\grcalcb) {\line(0,1){\dfactor}}
   \advance \grcalca by \hfactor
   \advance \grcalca by -\dfactor
   \advance \grcalcb by \dfactor
   \put(\grcalca,\grcalcb) {\line(-1,0){\factor}}
   \advance \grcalcb by \factor
   \advance \grcalcb by -\sfactor
   \put(\grcalca,\grcalcb) {\line(-1,0){\factor}}
   \grcalcc = \factor
   \advance \grcalcc by -\sfactor
   \put(\grcalca,\grcalcb) {\line(0,-1){\grcalcc}}
   \advance \grcolumn by 1}

 \newcommand{\grmpt}{
   \grcalca = \grcolumn
   \multiply \grcalca by \factor
   \advance \grcalca by \hfactor
   \grcalcb = \grrow
   \multiply \grcalcb by \factor
   \put(\grcalca,\grcalcb) {\line(0,-1){\dfactor}}
   \advance \grcalca by \hfactor
   \advance \grcalca by -\dfactor
   \advance \grcalcb by -\dfactor
   \put(\grcalca,\grcalcb) {\line(-1,0){\factor}}
   \advance \grcalcb by -\factor
   \advance \grcalcb by \sfactor
   \put(\grcalca,\grcalcb) {\line(-1,0){\factor}}
   \grcalcc = \factor
   \advance \grcalcc by -\sfactor
   \put(\grcalca,\grcalcb) {\line(0,1){\grcalcc}}
   \advance \grcolumn by 1}

 \newcommand{\grmpb}{
   \grcalca = \grcolumn
   \multiply \grcalca by \factor
   \advance \grcalca by \hfactor
   \grcalcb = \grrow
   \multiply \grcalcb by \factor
   \advance \grcalcb by -\factor
   \put(\grcalca,\grcalcb) {\line(0,1){\dfactor}}
   \advance \grcalca by \hfactor
   \advance \grcalca by -\dfactor
   \advance \grcalcb by \dfactor
   \put(\grcalca,\grcalcb) {\line(-1,0){\factor}}
   \advance \grcalcb by \factor
   \advance \grcalcb by -\sfactor
   \put(\grcalca,\grcalcb) {\line(-1,0){\factor}}
   \grcalcc = \factor
   \advance \grcalcc by -\sfactor
   \put(\grcalca,\grcalcb) {\line(0,-1){\grcalcc}}
   \advance \grcolumn by 1}

 \newcommand{\grmp}{
   \grcalca = \grcolumn
   \multiply \grcalca by \factor
   \advance \grcalca by \factor
   \advance \grcalca by -\dfactor
   \grcalcb = \grrow
   \multiply \grcalcb by \factor
   \advance \grcalcb by -\dfactor
   \put(\grcalca,\grcalcb) {\line(-1,0){\factor}}
   \advance \grcalcb by -\factor
   \advance \grcalcb by \sfactor
   \put(\grcalca,\grcalcb) {\line(-1,0){\factor}}
   \grcalcc = \factor
   \advance \grcalcc by -\sfactor
   \put(\grcalca,\grcalcb) {\line(0,1){\grcalcc}}
   \advance \grcolumn by 1}
\newcommand{\gsy}{
   \grcalca = \grcolumn
   \multiply \grcalca by \factor
   \advance \grcalca by \hfactor
   \grcalcb = \grcalca
   \advance \grcalcb by \hfactor
   \grcalcc = \grcalca
   \advance \grcalcc by \factor
   \grcalcd = \grrow
   \multiply \grcalcd by \factor
   \grcalce = \grcalcd
   \advance \grcalce by -\tfactor
   \grcalcf = \grcalcd
   \advance \grcalcf by -\hfactor
   \grcalcg = \grcalce
   \advance \grcalcg by -\tfactor
   \grcalch = \grcalcd
   \advance \grcalch by -\factor
   \qbezier(\grcalcc,\grcalcd)(\grcalcc,\grcalce)(\grcalcb,\grcalcf)
   \qbezier(\grcalcb,\grcalcf)(\grcalca,\grcalcg)(\grcalca,\grcalch)
   \advance \grcalcf by -\dfactor
   \advance \grcalcb by \sfactor
   \qbezier(\grcalcc,\grcalch)(\grcalcc,\grcalcg)(\grcalcb,\grcalcf)
   \qbezier(\grcalca,\grcalcd)(\grcalca,\grcalce)(\grcalcb,\grcalcf)
   \advance \grcolumn by 2}
 \newcommand{\gwmuh}[3]{
   \grcalca = \grcolumn
   \multiply \grcalca by \factor
   \grcalcb = #2
   \advance \grcalcb by #3
   \multiply \grcalcb by \qfactor
   \advance \grcalca by \grcalcb
   \grcalcb = \grrow
   \multiply \grcalcb by \factor
   \grcalcc = #3
   \advance \grcalcc by -#2
   \multiply \grcalcc by \hfactor
   \grcalcd = \factor
   \advance \grcalcd by \hfactor
   \put(\grcalca,\grcalcb){\oval(\grcalcc,\grcalcd)[b]}
   \grcalca = \grcolumn
   \multiply \grcalca by \factor
   \grcalcc = #1
   \multiply \grcalcc by \hfactor
   \advance \grcalca by \grcalcc
   \advance \grcalcb by -\hfactor
   \advance \grcalcb by -\qfactor
   \put(\grcalca,\grcalcb) {\line(0,-1){\qfactor}}
   \advance \grcolumn by #1}

 \newcommand{\gwcmh}[3]{
   \grcalca = \grcolumn
   \multiply \grcalca by \factor
   \grcalcb = #2
   \advance \grcalcb by #3
   \multiply \grcalcb by \qfactor
   \advance \grcalca by \grcalcb
   \grcalcb = \grrow
   \advance \grcalcb by -1
   \multiply \grcalcb by \factor
   \grcalcc = #3
   \advance \grcalcc by -#2
   \multiply \grcalcc by \hfactor
   \grcalcd = \factor
   \advance \grcalcd by \hfactor
   \put(\grcalca,\grcalcb){\oval(\grcalcc,\grcalcd)[t]}
   \grcalca = \grcolumn
   \multiply \grcalca by \factor
   \grcalcc = #1
   \multiply \grcalcc by \hfactor
   \advance \grcalca by \grcalcc
   \advance \grcalcb by \factor
   \put(\grcalca,\grcalcb) {\line(0,-1){\qfactor}}
   \advance \grcolumn by #1}

 \newcommand{\gsbox}[1]{
   \grcalca = \grcolumn
   \multiply \grcalca by \factor
   \grcalcb = \grrow
   \multiply \grcalcb by \factor
   \advance \grcalcb by -\factor
   \grcalcc = #1
   \multiply \grcalcc by \factor
   \grcalcd = \factor
   \put(\grcalca,\grcalcb){\framebox(\grcalcc,\grcalcd){}}}

   \newcommand{\gbox}[2]{
   \grcalca = \grcolumn
   \multiply \grcalca by \factor
   \grcalcb = \grrow
   \multiply \grcalcb by \factor
   \advance \grcalcb by -\factor
   \grcalcc = #1
   \multiply \grcalcc by \factor
   \grcalcd = #2
   \multiply \grcalcd by \factor
   \put(\grcalca,\grcalcb){\framebox(\grcalcc,\grcalcd){}}}

\allowdisplaybreaks[4]

\newcommand{\linea}{\gcl{1}}
\newcommand{\cruce}{\gbr}

\theoremstyle{remark}
\newtheorem*{remark}{Remark}

\font\cyrillic=wncyi10
\newcommand{\SU}{\mathop{\hbox{{\cyrillic UX}}}}

\title{Hopf algebras in non-associative Lie theory}

\author{J.~Mostovoy, J.M.~Perez-Izquierdo, I.P.~Shestakov}
\begin{document}

\maketitle

\begin{abstract} We review the developments in the Lie theory for non-associative products from 2000 to date and describe the current understanding of the subject in view of the recent works, many of which use non-associative Hopf algebras as the main tool.
\end{abstract}

\section{Introduction: Lie groups, Lie algebras and Hopf algebras}

Non-associative Lie theory, that is, Lie theory for non-associative products, appeared as a a subject of its own in the works of Malcev who constructed the tangent structures corresponding to Moufang loops. 
Its general development has been slow; nevertheless, by now many of the basic features of the theory have been understood, with much of the progress happening in the last ten years or so. In the present paper we outline the non-associative Lie theory in general and review the recent developments.

The history of the subject before 2000 has been summarized in the paper of Sabinin \cite{SabininHistory}. Probably, the single most important advance to that date has been the definition of a Sabinin algebra by Mikheev and Sabinin in 1986 (they use the term {\em hyperalgebra}). Sabinin algebras are the structures tangent to general non-associative unital local products; under some convergency conditions one can recover the product from a Sabinin algebra. Lie, Malcev, Bol, Lie-Yamaguti algebras and Lie triple systems are all particular instances of Sabinin algebras.

Classical Lie theory relates three kinds of objects: Lie groups, Lie algebras and cocommutative Hopf algebras. Until recently, however, Hopf algebras had little, if any, role in the non-associative Lie theory. To a certain extent, this gap has now been filled, and here we shall give an outline of the theory of local loops and Sabinin algebras from the point of view of non-associative Hopf algebras.  We shall be guided by the analogy with the theory of Lie groups and Lie algebras, so first we recall the basics of the classical theory.

At the heart of Lie theory lies the equivalence of the following three categories:
\begin{itemize}
\item simply-connected finite-dimensional Lie groups;
\item finite-dimensional Lie algebras;
\item irreducible cocommutative finitely generated Hopf algebras.
\end{itemize}
Let us describe explicitly the functors that establish this equivalence. In what follows all the vector spaces, Lie algebras etc will be defined over a field of characteristic zero.

\subsection{Lie groups $\to$ Lie algebras.}
There are several ways to produce the Lie algebra of a Lie group. Let $G$ be a Lie group with the unit element $e$. The most common construction identifies the tangent space ${T}_e G$ 
with the vector space of left-invariant vector fields on $G$. The commutator of two such fields (considered as derivations acting on functions) is again left-invariant, and defines the Lie bracket on ${T}_e G$.

This definition can be stated in different terms. On the tangent bundle ${T} G$ there exists a canonical flat connection $\nabla$. It can be defined by the associated parallel transport: for any path between $a,b\in G$ it sends ${T_a} G$ to ${T_b} G$ by the differential of the left multiplication by $ba^{-1}$. The curvature tensor of $\nabla$ is identically zero and the torsion tensor ${\rm T}$ is covariantly constant, that is $\nabla {\rm T}=0$ on $G$. This implies that the torsion tensor ${\rm T}$ is defined completely by its value on $T_e G$, which is a linear map $$T_e G\otimes T_e G\to T_e G.$$
This map coincides, up to sign, with the Lie bracket on $T_e G.$

\subsection{Lie algebras $\to$ Lie groups.}  One way to obtain a Lie group from a Lie algebra is via the Baker-Campbell-Hausdorff formula which expresses the formal power series
\begin{equation}\label{eq-bch}
\log (\exp x \exp y)
\end{equation}
in terms of iterated commutators in $x$ and $y$. For finite-dimensional Lie algebras, the Baker-Campbell-Hausdorff series always has a non-zero radius of convergence 
and, hence, defines a local Lie group.

It can be shown that finite-dimensional local Lie groups can always be extended to global Lie groups, but the proof of this fact is somewhat mysterious and uses the facts that seem to come from outside of Lie theory, such as the Ado theorem. The Ado theorem claims that every finite-dimensional Lie algebra is isomorphic to a Lie subalgebra of the algebra of the matrices. Therefore, it integrates to a local Lie group which is embedded as a local subgroup into the group of invertible matrices. Then it can be shown that this local subgroup generates a subset in the group of matrices which is an image of a smooth injective Lie group homomorphism. This last step is far from being trivial and uses in a crucial way the associativity of the matrix multiplication.

\subsection{Lie algebras $\to$ Hopf algebras.}
Each Lie algebra $\g$ can be embedded into an associative algebra $A$ in such a way that the bracket in $\g$ is induced by the commutator in $A$. Among all such embedding one is universal. The corresponding algebra (called universal enveloping algebra) is constructed as the quotient of the unital tensor algebra on $\g$ by the relations $$x\otimes y - y\otimes x - [x,y]$$
for all $x,y\in \g$ and is denoted by $U(\g)$. The Lie algebra $\g$ is naturally embedded into the tensor algebra on $\g$ and this embedding descends to the embedding
$$\g\to U(\g)$$
into the universal enveloping algebra.

The algebra $U(\g)$ can be given a coproduct $$\delta: U(\g)\to U(\g)\otimes U(\g)$$ which turns it into a Hopf algebra. It is defined by declaring all the elements of $\g\subset U(\g)$ to be primitive:
$$\delta(x)=x\otimes 1+ 1\otimes x.$$
Since the coproduct in a Hopf algebra is an algebra homomorphism, this is sufficient to determine $\delta$ completely. According to the Poincar\'e-Birkhoff-Witt theorem, as a coalgebra $U(\g)$ is isomorphic to the symmetric algebra $k[\g]$. This latter can be thought of as the universal enveloping algebra of the abelian Lie algebra which coincides with $\g$ as a vector space and whose Lie bracket is identically zero.

It is immediate that as a coalgebra $U(\g)$ is cocommutative. Also, by construction, it is generated by its primitive elements.

\subsection{Hopf algebras $\to$ Lie algebras.}

A straightforward calculation shows that the subspace of primitive elements in a Hopf algebra is always closed under the commutator.  In particular, there is a functor that assigns to each Hopf algebra the Lie algebra of its primitive elements. It is not hard to prove that it gives an equivalence of the category of cocommutative primitively generated Hopf algebras to that of Lie algebras (of not necessarily finite dimension).

\subsection{Lie groups $\to$ Hopf algebras.}
Let $D_e G$ be the space of distributions (that is, linear functionals on functions) on the Lie group $G$, supported at the unit $e$. In other words, $D_eG$ is the space spanned by the delta function supported at $e$ and all of its derivatives. The multiplication on $G$ gives rise to a convolution product on $D_e G$; the diagonal map $G\to G\times G$ induces a coproduct on $D_e G$. It can be seen that $D_eG$ is a Hopf algebra; moreover, it is precisely the universal enveloping algebra of the Lie algebra of $G$.

\subsection{Hopf algebras $\to$ Lie groups.} To recover a Lie group from a Hopf algebra, one uses the fact that any cocommutative Hopf algebra generated by its primitive elements is a universal enveloping algebra of some Lie algebra. In particular, by the Poincar\'e-Birkhoff-Witt theorem, it can be thought of as the symmetric algebra $k[V]$ on some vector space $V$ over the base field $k$ with a modified product $$\mu:k[V]\otimes k[V]\to k[V].$$
Consider the composition of this product with the projection $k[V]\to V$ onto the space of primitive elements. The resulting map $$F_\mu:k[V]\otimes k[V]\to V$$
can be interpreted as follows. Choose a basis  $e_1,\ldots, e_n$ in $V$. Then $k[V]$ is the algebra of the polynomials in the $e_i$. Each component of the map $F_\mu$ is a formal power series in $x_i$ and $y_j$, where 
$x_1,\ldots, x_n$ and $y_1,\ldots, y_n$ are the coordinates  dual to the $e_i$ in the first and the second copy of $V$ respectively. It turns out that for any Hopf algebra such that $V$ is of finite dimension these formal power series all converge in a neighbourhood $U$ of the origin in $V\times V$ and thus give an analytic map
$$V\times V\supseteq U\to V.$$
This map is, in fact, a finite-dimensional local Lie group and these are always locally equivalent to simply-connected Lie groups.

\subsection{The Ado theorem}

The fundamental triangle of Lie theory is remarkably robust: it can be generalized to a variety of situations. For instance, finite-dimensional Lie algebras can be replaced by infinite-dimensional Lie algebras, Lie algebras in tensor categories other than vector spaces, or, as we shall see, can be substituted for more general kinds of tangent algebras.  There is only one construction that does not always survive into the more general context.

Note that while recovering a Lie group either from its Lie algebra or from its Hopf algebra of distributions we first construct a formal Lie group, that is, a collection of power series which a priori may not converge. In order to see that this formal Lie group comes from an actual Lie group we need the Ado theorem, which states that every finite-dimensional Lie algebra has a faithful finite-dimensional representation. In particular, it follows that any finite-dimensional formal Lie group can be obtained from a local subgroup in a matrix group. A further argument is then needed to see that every such local subgroup extends to an injective homomorphism of a Lie group to the group of invertible matrices.

This proof breaks down already in the case of infinite-dimensional Lie groups. It also fails in the case of general finite-dimensional  loops. Nevertheless, there are situations (Moufang loops, nilpotent loops) where each formal multiplication extends to a global loop. In these situations we have an analogue of the Ado theorem, though we should point out that this theorem  alone does not guarantee the globalizability of formal products.

\subsection{Nilpotent groups}
Lie algebras and Hopf algebras also appear in the theory of discrete nilpotent groups.

Let $G$ be a discrete group. Write $\Q G$ for the group algebra of $G$ over the rationals and let $\Delta\subset \Q G$  be the augmentation ideal. The $i$th dimension subgroup $D_i G$ (or $D_i (G,\Q)$) consists of all those $g\in G$ for which $$g-1\in \Delta^i.$$
The dimension subgroups form a descending central series which is closely related to the lower central series $\gamma_i G$ defined inductively by $\gamma_1 G= G$ and $\gamma_i G= [\gamma_{i-1}G, G]$ for $i>1$. Namely, $D_i G$ consists of all $g\in G$ for which there exists $n\geq 1$ with the property that $g^n\in \gamma_i G$. This latter statement was proved by Jennings \cite{Jennings}.

The graded abelian group $$\mathcal{L}G= \bigoplus_{i>0}D_i G/D_{i+1} G$$ associated with the dimension series
is actually a Lie ring with the Lie bracket induced by the group commutator $$[a,b]=a^{-1}b^{-1}ab.$$ Tensored with the rational numbers, this Lie ring becomes a Lie algebra and its universal enveloping algebra can be explicitly identified as the algebra
$$\bigoplus_{i\geq 0} \Delta^i/\Delta^{i+1}$$
associated with the filtration on $\Q G$ by the powers of the augmentation ideal.

In general, a discrete group $G$ cannot be recovered from the Lie algebra $\mathcal{L} G\otimes \Q$ coming from the dimension series (equivalently, from the lower central series). However, if $G$ is a finitely generated nilpotent group,  $\mathcal{L} G\otimes \Q$ is finite-dimensional and can be integrated to a Lie group into which $G$ embeds as a discrete subgroup if it is torsion-free.

One group for which the dimension series and the lower central series can be described explicitly is the free group $F_n$ on $n$ generators $x_1,\ldots, x_n$. Let $\Z\langle\langle X_1,\ldots, X_n\rangle\rangle$ be the ring of formal power series in $n$ non-commuting variables. The group of invertible elements in this ring consists of power series which start with $\pm 1$. We denote this group by $\Z\langle\langle X_1,\ldots, X_n\rangle\rangle^*$.
There is a homomorphism
\begin{align}\label{eq:magnus}
\mathcal{M}: F_n & \to  \Z\langle\langle X_1,\ldots, X_n\rangle\rangle^*,\\
x_i&\mapsto  1+X_i,\nonumber
\end{align}
which is, in fact, injective. The $i$th terms of the dimension and of the lower central series of $F_n$ coincide and consist of those elements which map under $\mathcal{M}$ to the power series of the form
$$1+\text{terms of degree at least\ } m.$$
This, together with the injectivity of $\mathcal{M}$, shows that the free group $F_n$ is residually nilpotent.
\medskip

The non-associative versions of the above constructions and statements are the subject of the present paper. There are fundamental (and vast) parts of Lie theory which have not been understood yet in the non-associative context, most notably the representation theory. By no means we want to suggest that there are inherent obstructions to this, apart form time and effort.

\section{Local, infinitesimal and formal loops}

The non-associative Lie theory deals with the following three equivalent categories:
\begin{itemize}
\item formal loops;
\item Sabinin algebras;
\item irreducible cocommutative and coassociative non-associative Hopf algebras.
\end{itemize}

In contrast to the usual finite-dimensional Lie theory, there is no procedure to extend formal loops to global multiplications on manifolds. This, however, should not be seen as a problem since many non-associative products are inherently of local nature. The true analogue of a Lie group in the non-associative Lie theory is not a smooth loop, but, rather, a local loop, or a germ of a local loop. (We shall use the term ``infinitesimal loops'' for germs of local loops.) A formal loop produces a local loop whenever the power series that define it converge in some neighbourhood of the origin; verifying this is a problem from analysis and we shall not discuss it here.

\subsection{Local loops and smooth loops}
Let $M$ be a smooth finite-dimensional manifold\footnote{we assume that $\dim M>0$ is part of the definition of a smooth manifold.}. A {local multiplication} on $M$ (or, more precisely, on $U\subseteq M$)
is a smooth map
$$F: U\times U\to M,$$
where $U\subseteq M$ is a non-empty open subset. If there exists $e\in U$ with the property that
$$F|_{e\times U}=\mathrm{Id} (U)=F|_{U\times e},$$
the local multiplication $F$ is called {unital}, or a {local loop}. The point $e$ is referred to as the {unit}; if the unit exists, it is necessarily unique.
Often, when working with a fixed local loop $F$, we shall write $F(x,y)$ simply as $x\cdot y$ or $xy$.

For any local loop there exist two local multiplications $V\times V\to M$ with $V\subseteq U$, denoted by $x/y$ and $y\backslash x$. They are defined by
$$x/y \cdot y = x$$
and
$$y\cdot y\backslash x=x$$
for all $x,y\in V$, and, for obvious reasons, are called the {right} and the {left divison}, respectively. The existence of both divisions follows from the fact that the {right} and {left multiplication} maps
$$R_y=F|_{U\times y}:U\to M$$
and
$$L_y=F|_{y\times U}: U\to M$$
are close to the inclusion map $U\hookrightarrow M$ when $y$ is close to $e$. In particular, if $y$ is sufficiently close to $e$, both maps $R_y$ and $L_y$ are one-to-one and their images contain a neighbourhood of $e$. We take $V$ to be the largest neighbourhood on which both divisions are defined.

In general, there is no reason to expect that the right and the left divisions would be expressed as multiplications by an inverse of any kind.

A local loop is called a {smooth loop} if $U=V=M$ in the above definitions. Recall that a {discrete loop}, or, simply, a {loop}, is a set $M$ with a globally defined product for which there exists a unit element, and such that the left and the right multiplications by a fixed element of $M$ are bijections. In particular, for each loop there are globally defined right and left divisions. In these terms, a smooth loop can be defined as a loop which, at the same time, is a smooth manifold, with the multiplication and the divisions being smooth maps.

\subsection{Example: invertible elements in algebras}

Call an element $a$ of a unital algebra invertible if both equations $ax=1$ and $xa=1$ have a unique solution. The invertible elements of a finite-dimensional algebra over the real numbers form a local loop.

This local loop is not necessarily a loop. Consider, for instance, the Cayley-Dickson algebras $\mathbb{A}_n$ on $\R^{2^n}$. When $n>3$ there exist pairs of invertible elements in $\mathbb{A}_n$ whose product is not invertible.

\subsection{Example: homogeneous spaces}

Recall that a smooth manifold $M$ is a {homogeneous space} for a Lie group $G$ if $G$ acts transitively on $M$ and the action is smooth. Choose a point $e\in M$. Then we have a smooth map $$p: G\to M,$$ which sends $g\in G$ to $g(e)\in M$. This map identifies $M$, as a smooth manifold, with the left coset space $G/G_e$, where $G_e$ is the stabilizer of  $e$. Conversely, if $G$ is a Lie group and $H$ is a closed subgroup, the set $G/H$ of left cosets of $G$ by $H$ is a smooth manifold and a homogenous space for $G$.

Let $M$ be a homogeneous space for $G$ and $U\subseteq M$ a neighbourhood of a point $e\in M$. Assume that we are given a {section} of $p$ over $U$, that is, a smooth map $i: U\to G$ such that $i(e)$ is the unit in $G$ and $p \circ i = {\mathrm {Id}}(U).$  Then $M$ is a local loop, with the multiplication $U\times U\to M$ defined as
$$(x,y)\mapsto p\left(i(x)i(y)\right).$$
When $p$ is actually a homomorphism of Lie groups, that is, when $G_e$ is a normal subgroup in $G$, this local loop structure is the same thing as the product on $M$ restricted to $U\times U$.

There are many important examples of homogeneous spaces, among them spheres, hyperbolic spaces and Grassmannians.

\subsection{Infinitesimal loops}

Let $F$ and $F'$ be two local loops on $U\subseteq M$ and $U'\subseteq M'$ respectively. A {morphism} $F\to F'$ is a map $f:M\to M'$ such that $f(U)\subseteq U'$ and such that
$$F'(f(x), f(y))=f(F(x,y))$$ for
all $x,y\in U$.
With this notion of a morphism local loops form a category.

If a morphism $f$ from $F$ to $F'$ is an open embedding, we shall say that $F$ is a {restriction} of $F'$. In this situation, the local loops $F$ and $F'$ may be thought of as being ``locally the same'' near the unit. Certainly, they should not be distinguished by any infinitesimal algebraic structure at the unit, such as those arising in the classical Lie theory (think of formal groups or Lie algebras). With this in mind, we define {local equivalence}  as the smallest equivalence relation on the category of local loops under which open embeddings are equivalences. An equivalence class of local loops is called an {infinitesimal loop}. 

Strictly speaking, it is the infinitesimal loops and not the local loops that are the main subject in the non-associative Lie theory.
However, for the sake of simplicity, we shall speak of local loops rather than infinitesimal loops, wherever possible. This should not lead to confusion.

\subsection{Analytic local loops and formal loops}
Examples of local loops are very easy to construct via power series, since the only condition a multiplication of a local loop has to satisfy is the existence of the unit. Consider an $n$-tuple of power series $F_i(x_j, y_k)$ where $1\leq i,j,k\leq n$, and assume that all of them converge in some neighbourhood of the origin in $\R^{2n}$. The map
$$(x_1,\ldots, x_n, y_1,\ldots, y_n) \mapsto (F_1(x_j, y_k), \ldots, F_n(x_j, y_k))$$
defines a local loop on $\R^n$, with the origin as the unit, if and only if
\begin{equation}\label{1}
F_i(0,\ldots, 0, y_1,\ldots, y_n)=y_i
\end{equation}
and
\begin{equation}\label{2}
F_i(x_1,\ldots, x_n, 0,\ldots, 0)=x_i
\end{equation}
for all $i$.
These conditions mean that the coefficient of $x_j$ and $y_j$ in $F_i$ is equal to one if $i=j$ and vanishes otherwise. Moreover, a monomial of degree two or more can have a non-zero coefficient in $F_i$ only if it contains both the $x_j$ and $y_{j'}$ for some $j,j'$. The only further restrictions on the coefficients of the $F_i$ come from the requirement that the power series $F_i$ converge.

A local loop on an analytic manifold whose multiplication can be written in this form in some coordinate chart is called {analytic}. Similarly, an infinitesimal loop is analytic if it has an analytic representative.

The convergence of the power series  $F_i$ that specify an analytic loop is often irrelevant. For instance, a power series with no constant term can be substituted into another power series instead of a variable, and thus we do not need convergence in order to speak of the algebraic identities, such as associativity or Moufang identites, satisfied by an analytic loop.
In what follows we shall refer to an $n$-tuple of formal power series $F_i$ satisfying (\ref{1}) and (\ref{2}) as a {formal loop}.  

Let us state this definition in the form that does not make reference to any explicit basis in $\R^n$, nor to the fact that $n$ is finite.
For a vector space $V$ and a field $k$ write $k[V]$ for the symmetric algebra of $V$ over $k$. If a basis $e_\alpha$ is chosen for $V$, this is the commutative algebra of polynomials in the $e_\alpha$. Monomials of the form
\begin{equation}\label{3}
e_{\alpha_1}^{n_1}\ldots e_{\alpha_k}^{n_k}
\end{equation}
with ${\alpha_1}<\ldots<{\alpha_k}$ form a basis of $k[V]$. If $x_\alpha$ is the coordinate corresponding to $e_\alpha$, then any formal power series $F$ in the $x_\alpha$ can be thought of as a linear function on $k[V]$. This function sends the monomial (\ref{3}) to the coefficient of $x_{\alpha_1}^{n_1}\ldots x_{\alpha_k}^{n_k}$ in the series $F$.

Now, an analytic map from $V$ to a $k$-dimensional vector space $W$ is given by a $k$-tuple of power series in the $x_\alpha$, where $k=\dim W$. Given that each power series is just a linear function on $k[V]$, we can define a {formal map} from a vector space $V$ to a vector space $W$ as a linear map $k[W]\to V$ which annihilates the constants.

A {formal multiplication} is a formal map from $V\times V$ to $V$. Given that $k[V\times V]$ is canonically isomorphic to $k[V]\otimes k[V]$, this is the same a  linear map
$$k[V]\otimes k[V]\to V.$$
A formal multiplication $F$ is a {formal loop} if
$$F|_{1\otimes k[V]}=\pi_V=F|_{k[V]\otimes 1},$$
where $\pi_V:k[V]\to V$ be the projection of a polynomial onto its linear part.

\subsection{The canonical connection and geodesic loops} Let $F$ be a local loop on $U\subseteq M$ and let $V\subseteq U$ be a neighbourhood on which the left division is well-defined. The canonical connection $\nabla$ on the tangent bundle $T V$ has virtually  the same definition as for the Lie groups (apart from the fact that it is not defined on the whole of $M$): for $a,b\in G$ the corresponding parallel transport  sends ${T_a} G$ to ${T_b} G$ by the differential of the smooth map
$$x\mapsto b(a\backslash x).$$
It is clear that $\nabla$ is flat. There is no reason, however, for the torsion tensor of $\nabla$ to be constant.

It can be seen that each flat connection on a neighbourhood of a point comes from a local loop, called the geodesic loop of the connection. Indeed, a connection $\nabla$ on a manifold $M$ gives rise to the exponential map
$$\exp_a: T_a M\supseteq U_a \to M,$$
defined on a neighbourhood of the origin in $T_a M$, for each $a\in M$. If $U_a$ is chosen to be small enough, this map has an inverse, denoted by $\log_a$. Fix a point $e\in M$. Then there is a neighbourhood $U\subseteq M$ of $e$ such that
\begin{equation}\label{eq:monoaltproduct}
a\cdot b =\exp_a{(a\log_e{b})},
\end{equation}
is well defined for all $a,b\in U$. Here $a\log_e{b}$
stands for the parallel transport of the vector $\log_e{b}\in
{T}_e M$ to ${T}_a M$.
The operation $\cdot$ defines a local loop on $U\subseteq M$ and it is a straightforward check that $\nabla$ is its canonical connection. This loop is called the geodesic loop of $\nabla$. In order to be precise, we should consider the geodesic loop as an infinitesimal rather than local loop since there is a choice involved in the construction of the neighbourhood $U$.

It can be shown that the necessary condition for an infinitesimal loop (which is also sufficient in the case of analytic loops) to be a geodesic loop of a flat connection is that it should satisfy the right alternative property. Namely, it should be represented by a local loop such that
$$(a\cdot b)\cdot b= a\cdot (b\cdot b)$$
whenever both sides are defined. In particular, each local loop $F$ gives rise to a right alternative infinitesimal loop, which is the geodesic loop of the canonical connection of $F$.  One can recover a local loop $F$ from the corresponding geodesic loop with the help of the additional operation $\Phi$ defined by
\begin{equation}\label{eq:phi}
a\cdot b= a\times\Phi(a,b).
\end{equation}
Here $\cdot$ denotes the multiplication in $F$ and $\times$ is the product in the geodesic loop. $\Phi$ here can be any function of $a$ and $b$ that satisfies $\Phi(e,b)=b$ and  $\Phi(a,e)=e$.

The construction of the geodesic loop associated with a flat connection is due to Kikkawa \cite{Kikkawa} and Sabinin \cite{Sabinin_geodesic_loops}.

There is also a corresponding formal notion of a canonical connection.  Given a formal loop $F$ on a vector space $V$, the {
formal canonical connection} of $F$ is the restriction of $F$ to the subspace $$k[V]\otimes V\subset k[V]\otimes k[V].$$
Among all the formal loops that give rise to the same formal connection there is exactly one right alternative formal loop.

The notion of a geodesic loop is of crucial importance in non-associative Lie theory. Indeed, there are two constructions of the tangent structure to a local (infinitesimal, formal) loop and both are based on the fact that there is a flat connection associated with a local loop, see the next section. In fact, in the approach to the Lie theory taken by Sabinin in \cite{Sabinin_book} the role of non-associative analogues of Lie groups is played not by loops but by affinely connected manifolds. This involves additional algebraic structures such as odules; at the moment of the writing of the present paper the theory of odules has not noticeably advanced beyond the results of  \cite{Sabinin_book}.

\section{Sabinin algebras}

Many of the well-known generalizations of Lie algebras involve only one or two operations: Malcev algebras have one binary bracket,  Lie triple systems one ternary bracket, Bol and Lie-Yamaguti algebras one binary and one ternary bracket. In contrast, Sabinin algebras which are the most general form of the tangent structure for loops, have an infinite set of independent operations, and, as a consequence, admit an infinite number of axiomatic definitions. At the moment there are three different natural constructions for the set of operations in a Sabinin algebra. Two of these constructions are due to Sabinin and Mikheev and we shall review them in this section. The third set of operations, which appeared for the first time in \cite{ShU}  in the study of the primitive elements in non-associative bialgebras, will be considered in the next section. The complete set of axioms for this third set of operations is presently not known.

\begin{remark} The fact that the definition of Mikheev and Sabinin involved an infinite number of operations may have contributed to its slow acceptance. For some time Akivis algebras were considered as possible analogues of Lie algebras. The definition of an Akivis algebra involves two operations: an antisymmetric binary bracket and a ternary bracket, with only one identity that relates the two operations and generalizes the Jacobi identity. In spite of the elegance and simplicity of this definition, the category of Akivis algebras is not equivalent to that of formal loops and, hence, is not suitable as a basis for non-associative Lie theory. We shall consider the Akivis algebras in more detail in Section~\ref{section:akivis}.
\end{remark}

Both ways of deriving the structure of a Sabinin algebra from a non-associative product that were proposed by Mikheev and Sabinin consist of two steps. First, a local loop is replaced by the corresponding geodesic loop (which, in the analytic case is the same as a right alternative local loop), and then the tangent operations for a geodesic loop are extracted from its canonical connection. This results in two infinite sets of completely independent operations. Mikheev and Sabinin used the term ``hyperalgebra'' for the algebraic structure tangent to a geodesic loop and ``hyperalgebra with multioperators'' for the general tangent structure. In our terminology hyperalgebras with multioperators will be called Sabinin algebras, and hyperalgebras (which are the same as hyperalgebras with the trivial multioperator) will be called flat Sabinin algebras. The adjective ``flat'' here is supposed to reflect the fact that these are the tangent structures to general flat connections.

 \subsection{Flat Sabinin algebras from the torsion tensor of a connection}
Let $k$ be a field of characteristic zero. A flat Sabinin algebra is a vector space $V$ over $k$ together with a set of maps
$$V^{\otimes n+2}\to V$$
$$X_1\otimes\ldots\otimes X_n\otimes Y\otimes Z\mapsto \langle X_1, \ldots , X_n; Y, Z\rangle$$
for all integers $n\geq 0$, satisfying the following identities: 
$$\langle X_1, \ldots, X_n; Y,Z \rangle=-\langle X_1, \ldots, X_n; Z,Y \rangle,$$
\begin{multline*}
\langle X_1, \ldots, X_r,A,B, X_{r+1},\ldots, X_n; Y,Z \rangle-
\langle X_1, \ldots, X_r,B,A, X_{r+1},\ldots, X_n; Y,Z \rangle\\
+\sum_{k=0}^r\sum_{\alpha}\langle X_{\alpha_1}, \ldots,
X_{\alpha_k}, \langle X_{\alpha_{k+1}}, \ldots, X_{\alpha_r}; A,B
\rangle, X_{r+1},\ldots, X_n; Y,Z \rangle=0,
\end{multline*}
$$\sigma_{X,Y,Z}\left(\langle X_1,\ldots, X_r,X; Y,Z \rangle+
\sum_{k=0}^r\sum_{\alpha}\langle X_{\alpha_1}, \ldots, X_{\alpha_k};
\langle X_{\alpha_{k+1}}, \ldots, X_{\alpha_r}; Y,Z \rangle, X
\rangle\right)=0.
$$
Here $\alpha$ varies over the set of all bijections $\{1,\ldots r\}\to \{1,\ldots r \}$, $i\to\alpha_i$ such that
$\alpha_1<\alpha_2<\ldots<\alpha_k$, $\alpha_{k+1}<\ldots<\alpha_r$,
$k=0,1,\ldots, r$, $r\geq 0$, and $\sigma_{X,Y,Z}$ denotes the sum over all cyclic permutations of $X,Y,Z$.

\medskip

The above definition may seem complex, but it has the following geometric meaning. 

Consider an affine connection $\nabla$ at a point $e$ of a manifold $M$. It is characterized by two tensors: the curvature tensor and the torsion tensor which are related by the Bianchi identities. If $\nabla$ is flat (as it is in the case of a canonical connection) the curvature tensor is identically zero. If we assume that $\nabla$ is analytic, the torsion tensor $\rm T$ is uniquely determined in a neighbourhood of $e$ by its value and the values of all its covariant derivatives at this point. These derivatives are multilinear operations on $T_e M$. Set
$$\langle Y, Z \rangle = \mathrm{T}(Y,Z)$$
 and
$$\langle X_1,\ldots, X_n;  Y, Z \rangle = \nabla_{X_1}\ldots \nabla_{X_n}\mathrm{T}(Y,Z)$$ 
for all $n>0$. Then the brackets defined in this way are antisymmetric in $Y$ and $Z$, since the torsion tensor is antisymmetric, and satisfy the identities which come from the Bianchi identities and their derivatives, and from the fact that the curvature is zero. These identities are precisely those that appear in the definition of a flat Sabinin algebra.

From this construction it is clear that given a flat Sabinin algebra $\el$ one can reconstruct, locally, the corresponding flat connection, and, hence, the infinitesimal loop, provided that the operations of $\el$ define converging power series. The convergence conditions are explicitly stated by Mikheev and Sabinin in \cite{MS2, MS1}. In general, even if the convergence conditions are not satisfied, the geometric reasoning can be applied so as to produce a formal loop from any flat Sabinin algebra. 
 
 \subsection{Flat Sabinin algebras from pairs of Lie algebras}\label{section:lieenv}
 An alternative approach to flat Sabinin algebras is to define them as the algebraic structure that exists on a complement to a subalgebra in a Lie algebra.
 
Start with a direct sum decomposition of vector spaces 
\begin{equation}\label{eq:decomp}
\g=\h\oplus \el,
\end{equation}
 where $\g$ is a Lie algebra and $\h$ is a subalgebra. Let $\pi_\el$ be the linear projection map $\pi_\el:\g\to \el$ with $\h$ as the kernel. It will be convenient to introduce a simplified notation for the right-normed iterated bracket in $\g$:
$$[X_1, X_2, \ldots, X_n]=[X_1, [X_2,[ \ldots ,X_n]\ldots ]].$$
Define a sequence of multilinear operations\footnote{Our notation may be potentially confusing where it concerns the trilinear operation $(\cdot,\cdot,\cdot)$. Here this is not the associator.} on $\el$ by setting for each $n\geq 2$ and $X_1,\ldots, X_n\in \el$ 
$$(X_1,\ldots, X_n):=\pi_\el [X_1, \ldots ,X_n].$$ 
There are two kinds of relations satisfied by these operations. The relations of the first kind reflect the fact that they come from Lie brackets. Namely, the anti-symmetry of the Lie bracket
gives
\begin{equation}\label{eq:one}
(X_{1},X_2)+(X_2,X_{1})=0,
\end{equation}
and the Jacobi identity translates into the following identity:
\begin{equation}\label{eq:two}
(X_{1},X_{2},X_3)+(X_{2},X_{3},X_{1})+ (X_{3},X_{1},X_{2})=0.
\end{equation}
The identities of the second kind express the fact that $\h$ is a
subalgebra and not just a vector subspace. If $A,B\in \g$ are
arbitrary elements, then
$$\pi_\el [\pi_\el A -A, \pi_\el B-B]=0.$$
Setting $A=[X_1,\ldots,X_n]$ and $B=[Y_1,\ldots,Y_m]$  we get
\begin{multline*}
\pi_\el[[X_1,\ldots,X_n],[Y_1,\ldots, Y_m]]
+((X_1,\ldots,X_n),(Y_1,\ldots,Y_m) )\\
=((X_1,\ldots,X_n),\, Y_1,\ldots,Y_m)-((Y_1,\ldots,Y_m),\,
X_1,\ldots, X_n).
\end{multline*}
Now, the first term in this expression can be rewritten iteratively using the Jacobi identity in terms of the right-normed brackets so that the last relation takes the form
\begin{multline}\label{eq:three}
-((X_1,\ldots,X_n),(Y_1,\ldots,Y_m) )+((X_1,\ldots,X_n),\,
Y_1,\ldots,Y_m)-((Y_1,\ldots,Y_m),\,
X_1,\ldots,X_n)\\
=\sum_{\alpha}(-1)^{h(\alpha)} (X_{\alpha_1},\ldots X_{\alpha_n},
Y_1,\ldots, Y_m),
\end{multline}
where the summation is taken over all permutations $\alpha$ of the
set $\{1,\ldots,n\}$ for which there exists $s$ such that
$\alpha_i<\alpha_{i+1}$ for $i<s$, and $\alpha_i>\alpha_{i+1}$ for
$i\geq s$, and $h(\alpha)=n-s$. (Clearly, for any such $\alpha$ we
have that $\alpha_s=n$.)

We can now define a flat Sabinin algebra as a vector space with a sequence of multilinear operations 
$(X_1, \ldots, X_n)$ for all $n\geq 2$ satisfying the relations (\ref{eq:one})-(\ref{eq:three}). It can be shown that if $\el$ is a flat Sabinin algebra in this sense, there exists a decomposition (\ref{eq:decomp}) which gives rise to the operations on $\el$, \cite{MPS}.

The relation between flat Sabinin algebras in this sense and germs of flat affine connections is quite straightforward. Given such a connection $\nabla$ in the neighbourhood of a point $e$ on a manifold $M$, we can identify the tangent space $T_e M$ with the space of all $\nabla$-parallel vector fields. Set $\g$ to be the Lie subalgebra generated by the $\nabla$-parallel vector fields inside the Lie algebra of all vector fields, and $\h$ to be the subalgebra of $\g$ consisting of the vector fields that vanish at $e$. Then we have the direct sum decomposition
$$\g=\h\oplus T_e M,$$
and $T_e M$ acquires the structure of a flat Sabinin algebra.

One can, in fact, write down explicit formulae relating the two sets of operations in Sabinin algebras, see \cite{MS2, MS1}.
 
 \subsection{Multioperators}
Let $F$ be a right alternative local loop $F$ on $U\subseteq M$ and $$\Phi :U\times U\to M$$ a function such that $\Phi(e,b)=b$, $\Phi(a,e)=e$ and such that the Jacobian matrix  of $\Phi(a,b)$ with respect to $b$ is the identity when $b=e$.
Define  a new local loop ${F}_\Phi$ by setting $${F}_\Phi(a, b) = F(a, \Phi(a, b)).$$ The canonical connection of ${F}_\Phi$ coincides with that of $F$; in particular, $F$ is the geodesic loop associated with ${F}_\Phi$. (Strictly speaking, ${F}_\Phi$ is defined on a smaller neighbourhood than $U$. As usual, one can replace local loops here by infinitesimal loops.)

If $\Phi$ is analytic, one can write
\begin{equation}\label{eq-multi}
\Phi(a,b)=b+\sum_{m\geq 1, n\geq 2} \Phi_{m,n}(a,\ldots, a; b, \ldots, b),
\end{equation}
where $\Phi_{m,n}(X_1, \ldots, X_m; Y_1, \ldots, Y_n)$ is linear in each $X_j$ and $Y_j$ and invariant with respect to all the permutations of the $X_i$ and of the $Y_j$:
\begin{equation}\label{eq:multi}
\Phi_{m,n}(X_1, \ldots, X_m; Y_1, \ldots, Y_n)=\Phi_{m,n}(X_{\sigma(1)}, \ldots, X_{\sigma(m)}; Y_{\tau(1)}, \ldots, Y_{\tau(n)})
\end{equation}
 for all $\sigma\in \Sigma_m, \tau\in \Sigma_n$.
 
This motivates the following definition. A Sabinin algebra is a flat Sabinin algebra $\el$ together with a set of operations
$$\Phi_{m,n}: \el^{\otimes m+n}\to \el$$ 
 for $m\geq 1, n\geq 2$, satisfying the symmetry conditions~(\ref{eq:multi}). Given all that has been said about the flat Sabinin algebras, it is not hard to show that  Sabinin algebras form a category which is equivalent to the category of all formal loops.
 
\subsection{Free Sabinin algebras} The free Sabinin algebra on a set of generators  $S$ is the universal Sabinin algebra generated by $S$. It can be constructed as the vector space spanned by symbols corresponding to all possible operations 
$\langle  \cdot,\ldots, \cdot\, ; \cdot, \cdot\rangle$ and $\Phi_{m,n}(\cdot,\ldots, \cdot\, ;\, \cdot,\ldots,\cdot)$ and their compositions, whose arguments are elements of $S$, with the relations of a Sabinin algebra imposed. Another way to construct the free Sabinin algebras is by using the techniques described in the next section. Namely,  the free Sabinin algebra on a set of generators  $S$ is the set of primitive elements 
in the free non-associative algebra on the same set of generators. 

Several types of bases in free Lie algebras are known: Hall bases, Lyndon-Shirshov bases \cite{Reu}, right normed bases of Chibrikov \cite{ChiLie}. The construction of Lyndon-Shirshov bases was extended by Chibrikov \cite{Chi} to free Sabinin algebras. The dimensions of the graded summands of  a free Sabinin algebra on $n$ generators were calculated in \cite{Bre0}.

Free Sabinin algebras enjoy many other properties characteristic of free Lie algebras  \cite{Chi}. In particular, the variety of Sabinin algebras is Schreier (every subalgebra of a free Sabinin algebra is free). This implies that all automorphisms of finitely generated free Sabinin algebras are tame, and that the occurrence problem for free Sabinin algebras is decidable. Also, finitely generated subalgebras of free Sabinin algebras are residually finite and the word problem is decidable for the variety of Sabinin algebras.

 \subsection{Further remarks} 
 Just as in the associative case, there exists a Baker-Campbell-Hausdorff formula for Sabinin algebras. We shall discuss it in the next section.

The two sets of multilinear operations in a Sabinin algebra are completely independent and can be considered as Taylor expansions of two non-linear operations. Sabinin in \cite{Sabinin_book} takes this approach pointing out that it can be applied to non-analytic loops, and branding the constructions of \cite{MS2} and \cite{MS1} as obsolete. The book  \cite{Sabinin_book} was written before it was discovered that Sabinin algebras are important in the theory of non-associative bialgebras; in this latter context non-linear operations are not easy to deal with. 

\section{Non-associative Hopf algebras}

It would be fair to say that most of the constructions and results on Sabinin algebras mentioned in the previous section were the product of a quest for a Lie theory of flat connections. The subject of non-associative Hopf algebras has been developed in the context of the long-standing problem of finding an alternative enveloping algebra for an arbitrary Malcev algebra. While the existence of such an alternative envelope is still an open question, it turned out \cite{ShPI1} that Malcev algebras have universal enveloping algebras which are very similar in their properties to usual cocommutative Hopf algebras. Later, it turned out that a similar construction can be carried out for Bol algebras \cite{PI1} and, more generally, all Sabinin algebras \cite{PI2}.

Another motivation for the development of the machinery of non-associative Hopf algebras was the question of whether the commutator and the associator are the only primitive operations in a non-associative bialgebra. It appeared as a conjecture in \cite{HoStra}; if it were true, it would imply an important role for the Akivis algebras in non-associative Lie theory. This conjecture was answered in the negative in \cite{ShU}: in non-associative bialgebras, apart from the commutator and the associator, there exists an infinite series of independent primitive operations.  These operations, called the Shestakov-Umirbaev operations, are obtained from the associator by means of a procedure resembling linearization; they are closely related to the associator deviations in the nilpotency theory of loops (see Section~\ref{section:dev}).

The role of the non-associative Hopf algebras in the fundamental questions of Lie theory such as integration was clarified in \cite{MP3}. Hopf algebraic techniques are also relevant in other problems, such as the Ado theorem, which will be discussed in Section~\ref{section:ado}. 

\subsection{The definition and basic properties}
In what follows we shall consider algebras and coalgebras over a field $k$ of characteristic zero. 

A unital bialgebra is a unital algebra which, at the same time, is a counital coalgebra in such a way that the coproduct is  an algebra homomorphism (or, equivalently, such that the product is a coalgebra morphism). It will be convenient to use Sweedler's notation in which the coproduct is written as 
$$ \delta (x) =  \sum x_{(1)}\otimes x_{(2)}.$$
A unital bialgebra $H$ with the coproduct $\delta$ and the counit $\epsilon$ is called a (non-associative) Hopf algebra if is endowed with two additional bilinear operations, the right and the left division
$$/: H\times H\to H, \qquad\qquad \backslash: H\times H\to H,$$
$$(x,y)\mapsto x/y,\qquad\qquad\quad (x,y)\mapsto x\backslash y,$$
such that 
$$\sum (y x_{(1)})/ x_{(2)} = \epsilon(x) y = \sum (y / x_{(1)}) x_{(2)}$$
and
$$\sum  x_{(1)}\backslash( x_{(2)} y) = \epsilon(x) y = \sum  x_{(1)}( x_{(2)}\backslash y).$$
When the product in a Hopf algebra is associative and the coproduct is coassociative we get the usual notion of a Hopf algebra. In this case the antipode $S$ can be defined as $S(x)=1/x$ and, generally, we have 
 $$ x\backslash y = S(x)y\quad \text{and} \quad x/y = xS(y)$$
for all $x,y$.

The subspace of a Hopf algebra generated by the unit is a simple subcoalgebra. If this subspace is the only simple subcoalgebra, the Hopf algebra is called irreducible. In irreducible Hopf algebras the right and left division are uniquely determined. 

While non-associativity is an essential property of Hopf algebras that appear in Lie theory, the coproduct in many interesting cases is, actually, coassociative. For instance, the loop algebra $k L$ of a loop $L$ is always a coassociative Hopf algebra. In what follows by a Hopf algebra we shall understand a cocommutative and coassociative Hopf algebra, unless explicitly stated otherwise.

\subsection{Primitive elements and Shestakov-Umirbaev operations}
In any associative bialgebra the set of primitive elements is closed under the algebra commutator. This can be expressed by saying that the commutator is a primitive operation. Any other primitive operation can be written in terms of the commutators, in the following sense:  if an associative bialgebra $H$ is generated by a set $\{x_i\}$ where each $x_i$ is primitive, then each primitive element in $H$ is a linear combination of iterated commutators in the $x_i$. We express this by saying that the commutator forms a complete set of primitive operations.

In non-associative bialgebras there are other primitive operations, most notably the algebra associator
$$(a,b,c)=(ab)c-a(bc).$$
As discovered by Shestakov and Umirbaev in \cite{ShU}, the set of primitive operations obtained from commutators and associators is not complete. For instance, in the free non-associative algebra on one primitive generator $x$ the element
$$(x^2,x,x)-2x(x,x,x)=(x^2x)x-x^2x^2-2x(x^2x)+2x (xx^2)$$
is primitive, yet cannot be obtained from $x$ by taking linear combinations and compositions of commutators and associators. 

For $m,n\geq 1$ let $\underline{u}= x_1, \ldots, x_m $ and $\underline{v}=y_1, \ldots, y_n$ be sequences of primitive elements in a bialgebra and write 
$u=((x_1x_2)\cdots)x_m$ and $v=((y_1y_2)\cdots )y_n$ for the corresponding left-normed products. 
The {Shestakov-Umirbaev operations}  $p(x_1,\dots,x_m; y_1,\dots,y_n; z)$ are
defined inductively by
$$
(u,v,z) = \sum u_{(1)}v_{(1)}\cdot
p(\underline{u}_{(2)}; \underline{v}_{(2)}; z),
$$
where $(x,y,z)$ denotes the associator and $z$ is
primitive. 
Here Sweedler's notation is extended so as to mean that the sum is taken over all decompositions of the sequences  $\underline{u}$ and $\underline{v}$ into pairs of subsequences $\underline{u}_{(1)}, \underline{u}_{(2)}$ and $\underline{v}_{(1)}, \underline{v}_{(2)}$; the expressions $u_{(1)}$ and $v_{(1)}$ are the left-normed products of the elements of  $\underline{u}_{(1)}$ and $\underline{v}_{(1)}$, respectively.

For instance, the operation which corresponds to  $m=n=1$ is just the associator.  The operations corresponding to $m=2, n=1$ and $m=1, n=2$ are 
$$p(x_1, x_2; y; z)= (x_1 x_2, y, z) -x_1(x_2, y, z) - x_2(x_1, y, z)$$
and
$$p(x; y_1, y_2; z)= (x, y_1 y_2,  z) -y_1(x, y_2, z) - y_2(x, y_1, z),$$
respectively.

When the bialgebra in question is a Hopf algebra, the Shestakov-Umirbaev operations can be written with the help of the left division:
$$
p(x_1,\dots,x_m;y_1,\dots,y_n;z)= \sum
(u_{(1)}v_{(1)})\backslash (u_{(2)},v_{(2)},z).
$$
This form of the definition may be preferable since it expresses each operation directly via the associator.

In \cite{ShU} it is shown that all Shestakov-Umirbaev operations are primitive, and that together with the commutator they form a complete set of primitive operations.

In associative algebras the commutator satisfies the Jacobi identity; in particular, the primitive elements of a bialgebra $A$ form a Lie subalgebra of the commutator algebra of $A$. In the non-associative case, the Lie algebras are replaced by Sabinin algebras. Indeed, set
\begin{eqnarray*}
\langle 1;y,z\rangle &=& \langle y,z\rangle = -[y,z] = -yz + zy\\
\langle x_1,\dots,x_m;y,z\rangle &=& \langle
\underline{u};y,z\rangle =  -p(\underline{u},y,z)
+p(\underline{u},z,y)\\
\Phi^{SU}(x_1,\dots,x_m;y_1,\dots,y_n) &=&\\
&&\hskip -3cm  \frac{1}{m!}\frac{1}{n!}\sum_{\tau\in \Sigma_m,\sigma\in
\Sigma_n} p(x_{\tau(1)},\dots,x_{\tau(m)};y_{\sigma(1)},\dots,
y_{\sigma(n-1)};y_{\sigma(n)})
\end{eqnarray*}
with $u=((x_1x_2)\cdots)x_m$, $\Sigma_m$ the symmetric group on $m$
letters and $m\geq 1, n\geq 2$. Then the identities of a Sabinin algebra are satisfied for all $x_i, y_j, z$ in an arbitrary non-associative algebra, \cite{ShU}. 

In particular, we get a functor from non-associative algebras to
Sabinin algebras
$$
A\mapsto \SU(A)
$$
that assigns to an algebra $A$ the Sabinin algebra $\SU(A)$ on the same vector space.
In the case when $A$ is a bialgebra, the primitive elements of $A$ form a Sabinin subalgebra 
$$\text{Prim} (A)\subset \SU(A),$$ and we have a functor $$A\to \text{Prim} (A)$$
generalizing the corresponding functor from associative bialgebras (or Hopf algebras) to Lie algebras.

\subsection{Universal enveloping algebras for Sabinin algebras}

It was proved in \cite{PI2} that each Sabinin algebra $\el$ can be realized as the subspace of primitive elements in a certain non-associative irreducible Hopf algebra $U(\el)$. Moreover, the correspondence
$$\el\to U(\el)$$
is functorial. 

When $\el$ is a Lie algebra, the Hopf algebra $U(\el)$ is the usual universal enveloping algebra of $\el$. The properties of $U(\el)$ also closely mirror those of the associative enveloping algebras. There is a version of the Poincar\'e-Birkhoff-Witt theorem which says that $U(\el)$, as a coalgebra, is isomorphic to the symmetric algebra $k[\el]$. As a corollary, if $\el$ is given a basis $(e_i)$, the algebra $U(\el)$ is additively spanned by the left normed products $(\ldots (e_{i_1}e_{i_2})\ldots ) e_{i_n}$ with $i_k\leq i_{k+1}$. As a consequence, Hopf algebras $U(\el)$ are right and left Noetherian for finite-dimensional $\el$ and have no zero divisors.

Most importantly, we have the Milnor-Moore theorem which states that each irreducible Hopf algebra is the universal enveloping algebra of the Sabinin algebra of its primitive elements:
$$U=U(\text{Prim}(U)).$$
This implies that the universal enveloping algebra functor is an equivalence of the category of Sabinin algebras with that of irreducible Hopf algebras.

An explicit construction of $U(\el)$ can be found in \cite{PI2} or in Section~4.3 of \cite{MP3}. 

It seems plausible that for a flat Sabinin algebra $\el$ represented as complement to a Lie subaIgebra $\h$ in a Lie algebra $\g$ the algebra $U(\el)$ should have an explicit description in terms of $U(\g)$ and $U(\h)$. At the moment we only have a description of this kind for Malcev algebras and Lie triple systems (Section~\ref{section:associative}).

\subsection{Loops and Hopf algebras of distributions}

A distribution on a manifold $M$ is a linear functional on the functions on $M$. A distribution is said to be supported at a point if its value on a function only depends on the germ of the function at this point. Write $D_e M$ be the space of distributions on a manifold $M$ supported at the point $e\in M$. Under some reasonable conditions, which we shall assume to be satisfied, $D_e M$ as a vector space is spanned by the Dirac delta function at $e$ and all of its (higher) derivatives.

The space $D_e M$ has the natural structure of a coalgebra with the coproduct induced by the diagonal map $M\to M\times M$. As a coalgebra, it is isomorphic to the symmetric algebra $k[T_eM]$, and the space of primitive elements in $D_e M$ consists of all the derivatives along some vector in $T_eM$. In particular, it is cocommutative, coassociative and irreducible, since all these properties hold in $k[T_e M]$.

A local loop structure $U\times U\to M$ in a neighbourhood of $e$ induces a product
$$D_e U\otimes D_e U\to D_e M=D_e U.$$
This product is readily seen to be a coalgebra homomorphism, and, hence, $D_e U=D_eM$ is a Hopf algebra. If $L$ is a local loop on $U$, the notation $D_e L$ will stand for  $D_eU$ with this Hopf algebra structure.

In the case when a local loop $L$ is analytic, the product on $D_e L$ can be expressed in terms of the coefficients of the power series that define $L$. In particular, this makes it possible to define the Hopf algebra of distributions for any formal loop, and the explicit formulae in the formal case are quite straightforward.

Indeed, any formal map $$\theta: k[V]\to W$$
can be lifted to a unique coalgebra morphism
$$\theta':k[V]\to k[W],$$
which is defined as
$$\theta'(\mu) =
\sum_{n=0}^\infty \frac{1}{n!} \theta(\mu_{(1)})\cdots \theta(\mu_{(n)}) = \epsilon(\mu) 1 + \theta(\mu) + \cdots.$$ 
In particular, any formal loop
$$F:k[V\times V]=k[V]\otimes k[V]\to V$$
lifts to a product 
$$F': k[V]\otimes k[V]\to k[V],$$
which gives $k[V]$ the structure of a non-associative Hopf algebra. If the series that define $F$ converge, this Hopf algebra coincides with the Hopf algebra of the distributions $D_e F$.

By the Milnor-Moore and the Poincar\'e-Birkhoff-Witt theorems, each irreducible Hopf algebra can be considered as a product on $k[V]$ for some vector space $V$. Taking the primitive part of this product we get a formal loop:
$$k[V]\otimes k[V]\to k[V]\to \text{Prim}(k[V])=V.$$
This establishes an equivalence between the category of formal loops and that of irreducible Hopf algebras.

\subsection{The two structures of a Sabinin algebra on the tangent space to a loop}
Given a local loop $L$ on a manifold $M$, there are two ways to obtain the operations of a Sabinin algebra on $T_e M$. The first way was described by Sabinin and Mikheev: the brackets are the derivatives of the torsion tensor of the canonical connection of $L$ and the multioperator  $\Phi$ measures the failure of $L$ to be right alternative. The second possibility is to use the Sabinin algebra
 operations, as defined by Shestakov and Umirbaev, on the space of primitive elements in the Hopf algebra $D_e L$.

It turns out that these two Sabinin algebra structures are almost, though not quite, the same. Namely, the brackets in both Sabinin algebras coincide, while the Mikheev-Sabinin multioperator $\Phi$ is different from the Shestakov-Umirbaev multioperator $\Phi^{SU}$. At the moment it is not known how to express the Mikheev-Sabinin multioperator in terms of the Shestakov-Umirbaev primitive operations.

 \subsection{Exponentials, logarithms and the Baker-Campbell-Hausdorff formula}
In the associative case,  the exponential map is the unique, up to a rescaling, power series that sends the primitive elements of a complete Hopf algebra bijectively onto its group-like elements. Since the group-like elements in a Hopf algebra form a group, the Baker-Campbell-Hausdorff formula (\ref{eq-bch}) endowes the Lie algebra of the primitive elements with the structure of a group. This construction can be applied to any Lie algebra so as to obtain the corresponding formal Lie group.
 
Complete Hopf algebras can be defined in the non-associative context as well; the exponential and logarithmic power series in this situation are series in one non-associative variable. In contrast with the associative case, however, the exponential and the logarithm are not defined uniquely by the condition that they interchange the primitive and the group-like elements. In fact, this condition is satisfied by an infinite-dimensional affine family of power series all of which can be used to define a functor from Sabinin algebras to formal loops via a formula of  Baker-Campbell-Hausdorff type. For instance,  there exists a unique non-associative exponential series with the additional property that 
$$\exp{2x}=\exp{x}\exp{x}$$
and whose linear term is $x$. 
The first few terms of the corresponding Baker-Campbell-Hausdorff formula were calculated explicitly by Gerritzen and Holtkamp in \cite{GH}.They do not speak of Sabinin algebras and, in particular, do not discuss whether their formula provides any kind of integration.

A different exponential series was considered in \cite{MP3}. Replace in the usual exponential series the term $x^n$ with the left-normed product of $n$ copies of a non-associative variable $x$:
$$\exp x=1+x+\frac{1}{2!}x^2+\frac{1}{3!}x^2 x+ \frac{1}{4!}(x^2 x)x+\ldots$$
The exponential defined in this way is, essentially, the exponential map of the canonical connection on the loop of formal power series in one non-associative variable and non-zero constant term. 

The non-associative logarithm, which is defined as the inverse of the exponential, can be calculated explicitly and its coefficients involve Bernoulli numbers. Let $f:\mathbb{N}\to R$ be a sequence of elements of some ring $R$, indexed by the natural numbers. It can be extended to an $R$-valued function defined on all non-associative monomials in $x$. Given a  non-associative monomial $\tau$ define $f_{\tau}$ inductively as follows. For $\tau=x$ set $f_\tau=1$. If $\tau\neq x$, there is only
one way of writing $\tau$ as a product $(\ldots((x
\tau_1)\tau_2)\ldots)\tau_k$. Set
$$f_{\tau}=f_k\cdot f_{\tau_1}\cdot \ldots \cdot  f_{\tau_k}.$$
With this notation we have
$$\log(1+x)=\sum_{\tau} \frac{B_\tau}{\tau!} \cdot \tau, $$
where $\tau$ varies over the set of all non-associative monomials in $x$, and $B_\tau$ and $\tau!$ are the extensions of the Bernoulli numbers and the factorial, respectively, to the set of non-associative monomials \cite{MP3}.

Rather than use (\ref{eq-bch}), we define the formal power series $\mathcal{L}(x,y)$ in two non-associative variables $x$ and $y$ by the formula
$$\mathcal{L}(x,y)=\log\left(\sum_{m,n\geq 0} \frac{1}{m!n!}  x^m y^n  \right),$$
where all the products inside the logarithm are assumed to be left-normed. It can be seen that this series can be written as an infinite linear combination of primitive operations in $x$ and $y$. If applied to a Sabinin algebra $\el$, it gives a formal loop whose  tangent algebra has the same brackets as $\el$ but whose multioperator is zero. In particular, it provides formal integration for flat Sabinin algebras \cite{MPS}.

If $\el$ is not flat, one defines the Baker-Campbell-Hausdorff series to be
$$\mathcal{L}(x,\Phi(x,y)),$$
where $\Phi(x,y)$ is the power series defined by the multioperator as in (\ref{eq-multi}). This formula provides formal integration for arbitrary Sabinin algebras.

\subsection{Akivis algebras and their envelopes} \label{section:akivis}

An Akivis algebra is a vector space with one bilinear and one trilinear operation, denoted by $[\cdot, \cdot]$ and $(\cdot, \cdot,\cdot )$ 
respectively, such that
$$[x,y]=-[y,x]$$
and
$$[[x,y],z]+[[y,z],x]+[[z,x],y]=(x,y,z)+(y,z,x)+(z,x,y)-(x,z,y)-(z,y,x)-(y,x,z)$$
for all $x,y,z$. Any algebra is an Akivis algebra, with $[\cdot, \cdot]$ being the commutator and $(\cdot, \cdot,\cdot )$ the associator.  If $A$ is an algebra, denote by $Ak(A)$ the algebra $A$ considered as an Akivis algebra. Each Sabinin algebra is also automatically an Akivis algebra since the commutator and the associators are both primitive operations.

It was proved in \cite{Sh} that for each Akivis algebra $\mathfrak{a}$ there is a Hopf algebra $U_{Ak}(\mathfrak{a})$ together with an injective morphism of Akivis algebras $\mathfrak{a}\hookrightarrow U_{Ak}(\mathfrak{a})$ which is universal in the sense that any Akivis algebra morphism $\mathfrak{a}\to Ak(A)$ lifts to a unique algebra homomorphism $U_{Ak}(\mathfrak{a})\to A$.

If $\mathfrak{a}$ is actually a Sabinin algebra, $U_{Ak}(\mathfrak{a})$ is, in general, different from $U(\mathfrak{a})$. In particular, the image of  $\mathfrak{a}$ in $U_{Ak}(\mathfrak{a})$ is contained in the subspace of primitive elements, but not every primitive element of $U_{Ak}(\mathfrak{a})$ comes from  $\mathfrak{a}$.

\subsection{Dual Hopf algebras}

Inverting the arrows in the categorical definition of an associative Hopf algebra one arrives to the definition of the same structure. In the non-associative case this is no longer true.  

A unital bialgebra $H$ with the coproduct $\delta$ and the counit $\epsilon$ is called a dual Hopf algebra if is endowed with two additional linear operations $H\to H\otimes H$, called the right and the left codivision
$$x\mapsto  \sum x_{/1/}\otimes x_{/2/},\qquad\qquad\quad x\mapsto  \sum x_{\backslash 1\backslash}\otimes x_{\backslash 2\backslash},$$
such that 
\begin{equation*}
 \sum  x_{/1/(1)}\otimes  x_{/1/(2)} x_{/2/} = \sum x_{(1)}\otimes \epsilon(x_{(2)})  =  \sum x_{(1)/1/}\otimes x_{(1)/2/}x_{(2)}
\end{equation*}
and
$$ \sum x_{\backslash 1 \backslash}x_{\backslash 2 \backslash (1)}\otimes  x_{\backslash 2 \backslash  (2)} =  \sum \epsilon(x_{(1)})\otimes x_{(2)}  =   \sum x_{(1)} x_{(2)\backslash 1 \backslash}\otimes x_{(2)\backslash 2 \backslash}.$$

The dual of a finite-dimensional Hopf algebra is, rather tautologically, a dual Hopf algebra. In particular, the space of functions on a finite loop is a dual Hopf algebra, with the codivisions defined as
$$\sum f_{/1/} (x)\otimes f_{/2/}(y) = f(x/y)\qquad\text{and}\quad \sum f_{\backslash 1\backslash}(x)\otimes f_{\backslash 2\backslash}(y) = f(y\backslash x).$$
It is dual to the loop algebra, which is a cocommutative and and coassociative Hopf algebra.

A more general example is the space of functions on an algebraic loop.

An (affine) algebraic loop is an affine variety with the structure of a loop for which the multiplication and both divisions are morphisms (regular functions). Since the space $\mathcal{O}(X)$ of regular functions on any affine variety $X$ has the property that 
$$\mathcal{O}(X\times X) = \mathcal{O}(X)\otimes\mathcal{O}(X),$$ 
it follows that for an affine algebraic loop $X$ the space  $\mathcal{O}(X)$  is a commutative and associative dual Hopf algebra, with the codivisions defined in the same fashion as for finite loops. Moreover, over an algebraically closed field $k$ the category of affine algebraic loops is equivalent to that of finitely generated commutative and associative dual Hopf algebras. (The details are, essentially, the same as for affine algebraic groups.)

\begin{remark}
In principle, one may consider bialgebras with two divisions and two codivisions satisfying the idenitites that hold in Hopf algebras and in dual Hopf algebras. However, it can be shown that, in general, there is no way to introduce codivisions on the Hopf algebras that are of main interest here, namely, the universal enveloping algebras of Sabinin algebras.
\end{remark}

\section{Identities and special classes of loops}

\subsection{Identities in Hopf algebras coming from identities in loops}

Informally speaking, a {non-associative word} is an expression formed from several indeterminates by applying the multiplication and the right and the left divisions. (In particular, a monomial is a non-associative word formed using multiplication only.) If $w$ is  such a word, a local loop on $M$ is said to {satisfy the identity $w=e$} if for all the values of the indeterminates in $M$ for which $w$ is defined, the result of performing all the operations in $w$ is the unit. More generally, one can speak of the identities of the form $w_1=w_2$ where $w_1$ and $w_2$ are two words. As we already pointed out, identities also make sense in formal loops.

An identity satisfied by a local (formal) loop can be lifted to an identity in its algebra of distributions. If we start with a formal loop $F$ on a vector space $V$, then a non-associative word $\boldsymbol{w}$ in $n$ variables gives rise to a formal map $w:k[V]^{\otimes n}\to V$ which can be lifted to a coalgebra map $w': k[V]^{\otimes n}\to k[V]$. Then the identity 
\begin{equation}\label{id-loop}
\boldsymbol{w}_1=\boldsymbol{w}_2
\end{equation}
holds in $F$ if and only if we have 
\begin{equation}\label{id-distr}
w_1'=w_2'
\end{equation}
with any choice of arguments in $D_e F$. The identity (\ref{id-distr}) in the algebra of distributions is referred to as the linearization of the identity (\ref{id-loop}). 

The simplest example is the associativity in a formal loop which is necessary and sufficient for the associativity of its algebra of distributions. An analogous statement holds for the commutativity. A less trivial example is the left Moufang identity
$$x(y(x z))=((x y)x)z.$$ It translates into the identity
\begin{equation}\label{eq:linmoufang}
\sum
\mu_{(1)}(\nu(\mu_{(2)}\eta))= \sum ((\mu_{(1)}\nu)\mu_{(2)})\eta
\end{equation}
in the corresponding algebra of distributions. Non-associative Hopf algebras in which (\ref{eq:linmoufang}) is satisfied are called Moufang-Hopf algebras.

Geodesic loops are right alternative; their algebras of distributions satisfy
\begin{equation}\label{eq-ra}
\sum
\mu(\nu_{(1)}\nu_{(2)})= \sum
(\mu\nu_{(1)})\nu_{(2)}.
\end{equation}

\subsection{Identities in Sabinin algebras coming from identities in loops}

Identities in loops also produce identities in Sabinin algebras. In principle, there is a straightforward way to find all the identities in a Sabinin algebra that are implied by a given loop identity. Namely, the Baker-Campbell-Hausdorff formula expresses the loop product in terms of the operations in the corresponding Sabinin algebra; imposing an identity on a loop, therefore, imposes an infinite set of identities on its Sabinin algebra.  This set can be hugely redundant and choosing a subset consisting of independent identities is difficult in practice.
Nevertheless, this is not necessarily an impossible task: Kuzmin \cite{Kuzmin} used this method to show that a Malcev algebra integrates to a local Moufang loop. 

\medskip

Since loop identities translate into identities in Hopf algebras and in Sabinin algebras one can state various specializations of the general Lie theory for loops that satisfy a given set of identities. The most fundamental example of this is the Lie theory for right alternative loops which relates them to flat Sabinin algebras and Hopf algebras satisfying (\ref{eq-ra}).
There are many other such Lie theories for loops with identities; we describe some of them below.

\subsection{Moufang loops and Malcev algebras}

Historically, the first generalization of the Lie theory to non-associative products was the theory of Moufang loops and Malcev algebras. Recall that a Malcev algebra is a vector space with an antisymmetric bracket that satisfies
\begin{equation*}
[J(a,b,c),a] = J(a,b,[a,c])
\end{equation*}
where $$J(a,b,c) = [[a,b],c] + [[b,c],a] + [[c,a],b].$$
In this setup, we have the following three equivalent categories: simply-connected finite-dimensional Moufang loops, finite-dimensional Malcev algebras and finitely generated irreducible Moufang-Hopf algebras. The equivalence of first two categories was established in the works of Malcev \cite{Malcev}, Kuzmin \cite{Kuzmin} and Kerdman \cite{Kerdman}, while the third category was introduced only recently in \cite{ShPI1}.

An important feature of this case is the fact that finite-dimensional Malcev algebras always integrate to globally defined Moufang loops. While the original proof of this fact \cite{Kerdman} is rather inaccessible, the modern argument uses the theory of groups and algebras with triality and consists in translating the problem into the context of Lie groups and Lie algebras, see Section~\ref{section:associative} for more details. The Ado theorem also holds for finite-dimensional Malcev algebras \cite{ShPI1}, see more in Section~\ref{section:rep}. 

Moufang loops provide an example of how one can pass from loop identities to Sabinin algebra identities via the identities in the bialgebra of distributions. Given a Moufang loop $M$, the Moufang-Hopf identity (\ref{eq:linmoufang}), which holds in $D_e M$, implies that the tangent space $T_e M=\Prim {D_e M}$ belongs to the generalized alternative nucleus of  $D_e M$. For an algebra $A$  the generalized alternative nucleus is defined as 
$$
\Nalt(A) = \{ a \in A \mid (a,x,y) = - (x,a,y) = (x,y,a)\  \text{for all } {x,y\in A}\}.
$$
It is a Malcev algebra with the bracket given by  $[a,b] = ab - ba$. The rest of the Sabinin algebra operations in $T_e M$ can be expressed via the Malcev bracket  \cite{PI2}; as a consequence, the Malcev algebra structure is sufficient to reconstruct a local Moufang loop. 

A long-standing open problem proposed by Kuzmin \cite{Kuzmin}  is whether all Malcev algebras are special, that is, whether each Malcev algebra can be embedded into the commutator algebra of an alternative algebra. The formal integration procedure gives an embedding of a Malcev algebra into the generalized alternative nucleus of its universal enveloping algebra. However, the universal enveloping algebra of a Malcev algebra is not alternative, in general. One can define the universal alternative enveloping algebra of a Malcev algebra $\m$ as the maximal alternative quotient of $U(\m)$. Explicit calculations of both kinds of enveloping algebras for specific Malcev algebras can be found in \cite{Bre1, Bre2, Bre3}.

\subsection{Bol loops and Bol algebras}
The variety of left Bol loops is defined by the left Bol identity
$$x(y(xz))=(x (y x))z.$$
Its linearization is the identity
\begin{equation}\label{id-leftbolalg}
\sum \mu_{(1)}(\nu(\mu_{(2)}\eta))= \sum (\mu_{(1)}(\nu\mu_{(2)}))\eta.
\end{equation}
A Bol algebra is a vector space with one bilinear antisymmetric bracket and one trilinear bracket that satisfy the following identities:
\begin{eqnarray*}
&[a,a,b] = 0, & \\
&[a,b,c] + [b,c,a] + [c,a,b] = 0, &\\
&[x,y,[a,b,c]] =[[x,y,a],b,c] + [a,[x,y,b],c] + [a,b,[x,y,c]]
\end{eqnarray*}
and
\begin{equation*}
[a,b,[x,y]] = [[a,b,x],y] + [x,[a,b,y]] +[x,y,[a,b]] +
[[a,b],[x,y]].
\end{equation*}
The category of infinitesimal analytic finite-dimensional Bol loops is equivalent to the category of finite-dimensional left Bol algebras and to the category of irreducible Hopf algebras satisfying (\ref{id-leftbolalg}). While the correspondence between Bol loops and Bol algebras is well-studied (see \cite{MS2, NS, Sabinin_book}), the Hopf algebras in this context were first considered in \cite{PI1}.

Similarly to the case of Moufang loops, one can arrive to the Bol algebra identities via the distribution bialgebras. The identity (\ref{id-leftbolalg}) satisfied in the bialgebra of distributions shows that the primitive elements always belong to the left alternative nucleus of the distribution bialgebra, where the left  alternative nucleus of an algebra $A$ is defined as 
$$
\LNalt(A) = \{a \in A \mid (a,x,y) = - (x,a,y) \ \text{for all } {x,y\in A}\}.
$$
Any subspace of $\LNalt(A)$ which is closed under the operations $$[a,b] = ab - ba$$ and $$[a,b,c] = a(bc)-b(ac) - c[a,b]$$ is a (left) Bol algebra. This is the case, in particular, for  $T_e M$ inside $\LNalt(D_e M)$, where $M$ is a Bol loop  \cite{PI1}. All the other Sabinin algebra operations can be expressed in terms of these two brackets and, thus the structure of a Bol algebra is sufficient to reconstruct a local left Bol loop. 

The speciality problem for Malcev algebras has its generalization to Bol algebras. Namely, one may ask whether each Bol algebra is contained in a left alternative algebra as a Bol subalgebra. 
I.Hentzel and L.Peresi have shown that this is not the case \cite{HP}. Using computer-aided computations, they found an identity of degree eight which is satisfied in the Bol algebra of any left alternative algebra, but not in the free Bol algebra. 

An important property that holds for Malcev algebras but not for general Bol algebras is the Ado theorem: it was proved in  \cite{PI3} that it fails for  Lie triple systems (which are Bol algebras with the trivial binary bracket).

\subsection{Nilpotent loops and nilpotent Sabinin algebras}
Let us say that a loop $L$ is nilpotent of class $n$ if the $n+1$st term of the commutator-associator filtration (see Section~\ref{section:dev}) of $L$ is trivial: $$\gamma_{n+1}L=\{e\}.$$ This is not the standard definition of nilpotency in loops; however, there are strong indications that it is the correct one \cite{HVdL, M1, M2, MPS}. Since the commutator-associator filtration is defined in terms of words (namely, commutators, associators and associator deviations and their compositions) nilpotent loops of class $n$ form a variety. 

The corresponding Sabinin algebras, unsurprisingly, are the nilpotent Sabinin algebras of class $n$. Define a bracket of weight $k$ in a Sabinin algebra to be a multilinear operation in $k$ arguments formed by composing the brackets and the operations $\Phi_{m,n}$. Then a Sabinin algebra is nilpotent of class $n$ if all the brackets of weight greater than $n$ vanish in it. For a flat Sabinin algebra, nilpotency is the same as the nilpotency of the Lie algebra in which it is a summand according to the construction of Section~\ref{section:lieenv}, see \cite{MPS}.

The Lie theory of nilpotent loops and nilpotent Sabinin algebras exhibits various features typical of Lie and Malcev algebras. In particular, a nilpotent Sabinin algebra always integrates to a globally defined simply-connected nilpotent loop of the same class. This is due to the fact that the Baker-Campbell-Hausdorff series in this case consists of a finite number of terms and, hence, always converges. Also, the Ado theorem holds, both for Sabinin algebras and (globally) for simply-connected nilpotent loops (see Section~\ref{section:ado} for more details on the Ado theorem in the non-associative setting). While in the case of Malcev algebras one may think that these good properties are a consequence of Malcev algebras being ``close'' to Lie algebras, nilpotent Sabinin algebras may have an arbitrarily big (though finite) number of independent operations.

\subsection{Connections with quasigroups} In various contexts in geometry where non-associative algebra proves to be useful, it is quasigroups rather than loops that appear naturally. While we do not intend to discuss this subject in detail, there are two points that we should mention here.

One important application of quasigroups is the theory of symmetric spaces as developed by Loos \cite{L}; this is  where bialgebras  first appeared  within non-associative Lie theory. A symmetric space can be defined algebraically as a manifold $M$ with a product $\circ$ that satisfies, for all $x,y\in M$, the following identities:
 $$x\circ x = x,$$ $$x\circ (x\circ y) = y,$$ $$x\circ (y\circ z) = (x\circ y)\circ (x\circ z)$$ and such that the map $S_x\colon y \mapsto x\circ y$ has the unique fixed point $y=x$ in some neighbourhood of $x$. Geometrically, the product $x\circ y$ is the reflection of the point $y$ about the point $x$.
The fact that each $e\in M$ is idempotent is sufficient to define a bialgebra structure on $D_e M$ for each $e$ and linearize the identities in $M$ to obtain identities in $D_e M$. For instance, if $M$ is a symmetric space, in each $D_e M$  we have:
\begin{equation}\label{id-lts}
\sum \mu_{(1)}\mu_{(2)} = \mu, \qquad \sum \mu_{(1)}(\mu_{(2)}\nu) = \epsilon(\mu)\nu, \qquad \mu(\eta\nu) = \sum (\mu_{(1)}\eta)(\mu_{(2)}\nu).
\end{equation}
The local uniqueness of the fixed point for the multiplication map $S_e$ implies that the right multiplication in  $D_e M$ by the only group-like element 
$$e \colon f \mapsto f(e)$$
is bijective. This, together with the identities (\ref{id-lts}), can be used to deduce algebraically that the tangent space  $T_e M = \Prim (D_e M)$ is a Lie triple system with the product 
$[a,b,c] = a(bc)-b(ac)$.

\medskip

The second point is related to the terminology concerning the identities in loops. Given a local quasigroup, that is, a locally invertible map $M\times M\to M$ defined in a neighbourhood of a point $(a,b)$, one can define a local loop on a neighbourhood of $a\cdot b\in M$ by 
$$xy= (x/b )\cdot(a\backslash y),$$
where $\cdot, \backslash$ and $/$ are the product and the two divisions in the quasigroup. This construction associates with a local quasigroup a family of local loops, and it may be useful to compare the identities that hold in each local loop. It turns out that the Moufang identity is satisfied at the same time in each of its local loops; in particular, it is sufficient to verify this identity in one of the local loops in order to conclude that it holds in all of them. In contrast, the monoassociativity identity
\begin{equation}\label{id-monoass}
(xx)x=x(xx)
\end{equation}
may be satisfied in one of the local loops but not in the others. A monoassociative quasigroup is defined as a local quasigroup all of whose local loops satisfy (\ref{id-monoass}). (This definition is of importance in web theory where monoassociative quasigroups correspond to hexagonal three-webs.) As a consequence, a monoassociative local loop is sometimes  defined as a local loop of a monoassociative quasigroup \cite{AkG, BreMa}. This is a much stronger condition on a loop than just the identity (\ref{id-monoass}). In particular, each right alternative local loop  satisfies (\ref{id-monoass}); note that flat Sabinin algebras have an infinite number of independent operations. On the other hand, Sabinin algebras tangent to monoassociative loops have three independent operations: one bilinear, one trilinear and one quadrilinear \cite{Mi96, S89b}. The complete set of relations for these operations is not known.

\section{Constructions involving associative Hopf algebras}\label{section:associative}

Since each Sabinin algebra can be thought of a subspace of  a Lie algebra, it is not surprising that non-associative Hopf algebras can be interpreted in terms of associative Hopf algebras. The construction of the universal enveloping algebra of a Sabinin algebra given in \cite{PI2} is based on such an interpretation. In general, this construction does not produce closed formulae; however, there are important special cases when it can be significantly simplified. This happens, for instance, for Moufang-Hopf algebras and for universal enveloping algebras of Bol algebras.

\subsection{Groups with triality}
One important step in the study of Malcev algebras was to establish its relation with groups with triality \cite{Glau,Doro,Miheev}. Each self-map $d$ of the octonion algebra, antisymmetric with respect to its standard quadratic form, gives rise to two antisymmetric maps $d',d''$ uniquely defined by the relation 
$$d(xy) = d'(x)y + x d''(y).$$
This phenomenon is known as the local triality principle and, in general, it is often related in one or another way with exceptional behaviour such as that of exceptional Lie algebras or Jordan algebras \cite{Jacobson}. Behind it, there is a certian representation of the symmetric group $\Sigma_3$ on three letters, acting on the Dynkin diagram of the multiplication Lie algebra of the split octonions, which is a central simple Lie algebra of type $D_4$ \cite{Scha}. This symmetry lifts to the corresponding group $\Spin(8)$ and gives it the structure of a group with triality. More generally, an abstract group with triality is a group $G$ together with two automorphisms $\sigma, \rho$ such that $$\sigma^2 = \rho^3 = \mathrm{id},\quad \sigma\rho\sigma = \rho^2,$$ and
$$ P(g)\rho(P(g))\rho^2(P(g)) = 1,$$
where $P(g) = g\sigma(g)^{-1}$. The automorphisms $\sigma$ and $\rho$ define a representation of $\Sigma_3$, which is independent on a particular choice of the generators $\sigma$ y $\rho$ of orders $2$ and $3$, respectively. If $H\subseteq G$ is the subgroup fixed by $\sigma$, the quotient $G/H$ is a Moufang loop with the product 
$$g_1 H * g_2 H = \rho\sigma(P(g_1))g_2 H.$$

Let $\g$ and $\h$ be the Lie algebras of $G$ and $H$, respectively; the group $\Sigma_3$ acts on $\g$ with $\h$ being fixed by $\sigma$. The Malcev algebra $\m$  tangent to $G/H$ is the complement to $\h$ in $\g$ consisting of the eigenvectors of $\sigma$ with eigenvalue $-1$. The bracket in $\m$ is defined by 
$$[x,y]_\m =  [\rho^2(x)-\rho(x),y].$$
This is the construction of a Lie algebra with triality. Each Malcev algebra arises as the $-1$-eigenspace of $\sigma$ in such a Lie algebra.

In order to construct the universal enveloping algebra $U(\m)$, embed $G/H$ into $G$ by means of the map
$$gH \mapsto \rho(\sigma(P(g))),$$ 
which is a global section of the principal fibration  $G\to G/H$. The image of this embedding is the set 
$$\mathcal{M}(G)=\{ P(g) \mid g \in G\}.$$
The Moufang product on $\mathcal{M}(G)$ which comes from $G/H$ can be expressed directly in terms of the product on $G$ as
$$x*y = \rho(x)^{-1}y\rho^2(x)^{-1},$$ 
see \cite{GrZ}. 
Linearizing this formula we get that the subcoalgebra 
$$\mathcal{M}(U(\g)) = \left\{ \sum x_{(1)}\sigma(S(x_{(2)}))\mid x \in U(\g)\right\},$$ where $S$ is the antipode of $U(\g)$, with the product
\begin{displaymath}
x*y = \sum \rho(S(x_{(1)})) y \rho^2(S(x_{(2)}))
\end{displaymath}
is isomorphic to $U(\m)$, see \cite{BeMaP}. 

One can speak of the categories of Lie groups with triality, Lie algebras with triality and Hopf algebras with triality. The constructions of this paragraph establish their equivalence to the categories of analytic Moufang loops, Malcev algebras and Moufang-Hopf algebras, respectively. In this way, for instance, the global integrability of finite-dimensional Lie algebras implies global integrability of finite-dimensional Malcev algebras.

\subsection{Bruck loops}

A similar construction works for Bruck loops. These are left Bol loops that satisfy the automorphic inverse identity
$$(xy)^{-1} = x^{-1}y^{-1},$$ 
where $x^{-1}$ is shorthand for $e/x$. The binary operation in the corresponding Bol algebra in this case is trivial, and the structure of a Bol algebra reduces to that of a Lie triple system. This shows that Bruck loops are related to symmetric spaces. Indeed, the product 
$$x\circ y = x(y^{-1}x)$$
endowes an arbitrary local Bruck loop with the structure of a locally symmetric space. Conversely, given a point $e$ in a symmetric space the Bruck loop product is recovered by means of the relation
$$xy \circ e = x \circ (e \circ( y\circ e)).$$
While the product of a locally symmetric space can be globalized, this is not true for the corresponding local loop since the map $x\mapsto x\circ e$ may be only locally, and not globally, bijective. 

A symmetric space can be considered as a quotient  $G/H$, where $G$ is a Lie group with an involutive automorphism $\sigma$ and $H$ is an open subgroup in the fixed set $G^{\sigma}$. Locally, the quotient $G/H$ can be embedded into $G$ with the help of the map
$$gH \mapsto P(g) = g\sigma(g)^{-1},$$ 
which is a local section of the quotient map $G\to G/H$. The image of this map coincides, locally, with the set
$$\mathcal{M}(G) = \{P(g)\mid g \in G\}.$$ It can be shown that in a neighbourhood of the unit in  $\mathcal{M}(G)$ each element has a unique square root. The product of a local Bruck loop on  $\mathcal{M}(G)$ can be written then 
as $$x*y = \sqrt{x}\, y\, \sqrt{x}$$ for $x,y \in \mathcal{M}(G)$. 

Similarly, a Lie triple system can be considered as the $-1$-eigenspace of an involutive automorphism in a Lie algebra $\g$. As a consequence, the universal enveloping algebra of a Lie triple system, viewed as a Sabinin algebra, can be defined in terms of the Hopf algebra  $U(\g)$ and the automorphism $\sigma$ as the coalgebra $$\mathcal{M}(U(\g)) = \left\{ \sum x_{(1)}\sigma(S(x_{(2)})) \mid x \in U(\g)\right\},$$
with the product
$$x*y = \sum r(x_{(1)})y r(x_{(2)})$$
where $r$ is the map corresponding to the linearization of the square root, see  \cite{MP2}.

\medskip

Moufang loops and Bruck loops are only two examples of the situation when the algebra of distributions can be constructed with the help of an associative Hopf algebra. In general, we get such a construction whenever we have a lifting of a local product in $G/H$ to the product in $G$ by means of a local section.

\section{Representation theory}\label{section:rep}
Compared to the representation theory of Lie groups and Lie algebras, the representation theory of loops and Sabinin algebras is still in its infancy. It is, probably, too early to give a coherent picture of the field; here we touch only on a few topics.

The theory of loop representations has been put on firm foundations by J.D.H.~Smith \cite{Smith}; the same approach works for Sabinin algebras and non-associative Hopf algebras. It rests on the observation of J.~Beck \cite{Beck} that the concept of a split-null extension has a natural definition in the categorical setting: if $X$ is an object in a category  $\mathcal{C}$, an $X$-module is an abelian group in the comma category  $\mathcal{C}\downarrow X$. The objects of this category are morphisms  $Y \stackrel{\pi}{\rightarrow} X$ and the morphisms are commutative diagrams of the form

\begin{center}
\begin{tikzcd}[column sep=small]
Y\arrow{dr}[swap]{\pi}\arrow{rr}&   &Z\arrow{dl}{\pi} \\
& X &
\end{tikzcd}
\end{center}

This categorical description encompasses all reasonable definitions of modules; in particular, that of a group representation and that of a module over a Lie algebra. It has an advantage of producing immediate bijections between the representations of a formal group, of its Lie algebra and of the universal enveleping algebra of the latter. The same is true in the non-associative context: the equivalence of the category of formal loops to that of Sabinin algebras and to the category of irreducible Hopf algebras identifies that the corresponding sets of representations. Let us spell out what Beck's definition means in practice.

\subsection{Loop representations}
Let  $\mathbf{Loops}$ be the category of all loops and let $L$ be an arbitrary loop. Introduce the following notation:
\begin{itemize}
\item $L[X]$ --- the free product of $L$ with the free loop on one generator $X$;
\item $L_a$ and $R_a$ --- the operators of left and right multiplication, respectively, by $a\in L$ in $L[X]$; 
\item $\U(L)$ --- the group of self-maps of $L[X]$ generated by the $L_a$ and $R_a$ for all $a\in L$; 
\item $\U(L)_e$ --- the stabilizer of the identity in $L[X]$. 
\end{itemize}
 Abelian groups in the category $\mathbf{Loops}\downarrow L$ are in one-to-one correspondence with the representations of the group $\U(L)_e$. Namely, given a $\U(L)_e$-module $E_e$, define
 $$E=E_e\times L$$
as a set, endowed with the product 
\begin{equation}\label{product}
(x,a)(y,b) = (r_{a,b}(x) + s_{a,b}(y) , ab), 
\end{equation}
where $$r_{a,b} = R^{-1}_{ab}R_bR_a\text{\ and\ } s_{a,b} = R^{-1}_{ab}L_aR_b.$$ 
Then the projection of the loop $E$ onto $L$ is an abelian group in $\mathbf{Loops}\downarrow L$. The sum, inverse and the neutral element that define on $E\rightarrow L$ the structure of an abelian group are induced by the corresponding operations in $E_e$. Each abelian group in $\mathbf{Loops}\downarrow L$ arises in this way, with the module $E_e$ defined uniquely up to an isomorphism.

If $L$ belongs to a variety of loops $\mathbf{V}$, the category $\mathbf{Loops}$ in this construction can be replaced by $\mathbf{V}$. The free loop and the free product in this case should be taken in $\mathbf{V}$; instead of the group $\U(L)_e$ one obtains a group $\U(L;\mathbf{V})_e$. The abelian groups in  $\mathbf{V}\downarrow L$ are then in correspondence with certain representations of $\U(L;\mathbf{V})_e$. For instance, if $\mathbf{V}=\mathbf{Groups}$ is the variety of groups, we have 
$$\U(L;\mathbf{Groups})_e \simeq L.$$
Given a representation 
$$\rho:L\to \Aut (V)$$
the product (\ref{product})
takes the form
$$(x,a)(y,b) = (x + \rho_a(y) , ab).$$ 
See \cite{Smith} for further details.

\subsection{Hopf and Sabinin algebra representations}
Let $\mathcal{H}$ be the category of all  irreducible Hopf algebras and $H$ an object in $\mathcal{H}$. Beck's definition translates to this context in the manner similar to the case of loops. 

An abelian group $A \stackrel{\pi}{\longrightarrow} H$ in the category $\mathcal{H}\downarrow H$ is a morphism such that the Hopf algebra $A$  factorizes as $$A \cong k[V] \otimes H,$$ where $k[V]$ is the symmetric algebra of a vector space $V$ and  the product in $A$ is given by 
$$(p\otimes a)(q\otimes b) = \sum (r_{a_{(1)},b_{(1)}}(p))(s_{a_{(2)},b_{(2)}}(q))\otimes a_{(3)}b_{(3)}$$
in terms of the action of some operators $r_{a,b}, s_{a,b}$ on $V$ whose action is extended to $k[V]$ as $$r_{a,b}(pq) = \sum r_{a_{(1)},b_{(1)}}(p)\, r_{a_{(2)},b_{(2)}}(q)$$
and, similarly,
$$s_{a,b}(pq) = \sum s_{a_{(1)},b_{(1)}}(p)\, s_{a_{(2)},b_{(2)}}(q).$$
The subalgebra $k[V]$ appears naturally as the equalizer of the coalgebra morphisms $\pi\colon A \rightarrow H$ and 
\begin{align*}
\epsilon 1\colon A& \rightarrow H\\
x&\mapsto \epsilon(x)1.
\end{align*}
The sum, inverse and the neutral element in $A\stackrel{\pi}{\longrightarrow} H$ are given by
\begin{center}
\begin{tikzcd}[column sep=small]
\sum (p\otimes a_{(1)})\otimes (q\otimes a_{(2)})\arrow{dr}[swap]{\pi}\arrow{rr}{+}&   &pq\otimes a\arrow{dl}{\pi} \\
& \epsilon(p)\epsilon(q)a &
\end{tikzcd}
\begin{tikzcd}[column sep=small]
(p,a)\arrow{dr}[swap]{\pi}\arrow{rr}{-}&   &(S(p),a)\arrow{dl}{\pi} \\
& \epsilon(p)a &
\end{tikzcd}
\begin{tikzcd}[column sep=small]
a\arrow{dr}[swap]{\pi}\arrow{rr}{0}&   &1\otimes a\arrow{dl}{\pi} \\
& a &
\end{tikzcd}
\end{center}
where $S$ is the antipode in $k[V]$. The category $\mathcal{H}\downarrow H$ has finite products:
The product $A_1 \otimes_H A_2$ of $A_1 \stackrel{\pi_1}{\longrightarrow} H$ and $A_2 \stackrel{\pi_2}{\longrightarrow} H$ is the equalizer of the morphisms  $\pi_1 \otimes \epsilon$ and $\epsilon \otimes \pi_2$.

Usually, in the associative context it is the vector space $V$ that is called an $H$-module and the abelian group $A\stackrel{\pi}{\longrightarrow} H$ is not given a special name. In this case the action of $r_{a,b}$ is trivial, that is, $r_{a,b}(p) = \epsilon(a)\epsilon(b)p$, while $s_{a,b}(q) = \epsilon(b)s_{a,1}(q)$. If we write $a\cdot q = s_{a,1}(q)$, the formula for the product in $A$ becomes 
$$(p\otimes a)(q\otimes b) = \sum p(a_{(1)}\cdot q) \otimes a_{(2)}b,$$
which is nothing but the smash product of the $H$-module algebra  $k[V]$ with $H$, where the $H$-module algebra structure on $k[V]$ is induced by the $H$-module structure on $V$. 

\medskip

Modules over Sabinin algebras are somewhat easier to describe than modules over loops and Hopf algebras. Beck's notion of a module for a Sabinin algebra $\el$ coincides with the usual idea of a split-null extension. In other words, a module over $\el$ is a Sabinin algebra $V \oplus \el$ such that $\el$ is a subalgebra and such that each multilinear operation vanishes if at least two of its arguments belong to $V$.

\subsection{An extended concept of a representation}

Sometimes there is rationale for enlarging the class of the extensions considered in the definition of the modules over a loop (Hopf algebra, Sabinin algebra). For instance, the requirement that the extension $\pi:E{\rightarrow} L$ lies within the variety $\mathcal{M}$ of Moufang loops leads to conditions on $r_{a,b}$ and $ s_{a,b}\in \U(L;\mathcal{M})_e$; these conditions are expressed as the vanishing of certain elements of the group algebra of $\U(L;\mathcal{M})_e$ while acting on $\pi^{-1}(e)$. These elements may fail to vanish on a tensor product of $\U(L;\mathcal{M})_e$-modules. As a consequence, in contrast to the case of groups, a tensor product of $L$-modules in general will not be a $L$-module. In the same fashion, a tensor product of Sabinin algebra representations cannot be expected to be a representation. 

The failure of the tensor product of two representations to be a representation is a basic reflection of non-associativity. In case of the Moufang loops this problem can be overcome by studying abelian groups $E\to L$ such that $E$ is not necessarily a Moufang loop, but requiring only that the image of the zero section $L \to E$ is contained in the set of Moufang elements of $E$. These are the elements $a$ such that $$a(x(ay)) = ((ax)a))y\text{\quad  and\quad } ((xa)y)a = x(a(ya))$$ for all $x,y \in E$; Moufang elements in any loop form a Moufang loop. In the context of Hopf algebras this corresponds to studying abelian groups $A\stackrel{\pi}{\longrightarrow}H $ such that in the decomposition $A \cong k[V] \otimes H$ the subspace $\Prim(H)$ (but not $V$) is contained in the generalized alternative nucleus of $A$. As for Malcev algebras, representations of this kind correspond to split-null extensions $(V \oplus \m,[\,,\,])$ of the Malcev algebra $\m$ such that the bracket $[\,,\,]$ is anticommutative and the adjoint map $$r_a \colon x \mapsto r_a(x) = [x,a]$$ satisfies
\begin{equation}
\label{eq:rep_Malcev}
[[r_a,r_b],r_c] = - [r_{[a,b]},r_c] + r_{[[a,b],c]+[[a,c],b]+[a,[b,c]]}.
\end{equation}
Not only the resulting representation theory behaves well with respect to tensor products but also among its other good properties is the existence of faithful finite-dimensional representations for finite-dimensional Malcev algebras, see \cite{ShPI1}.

In general, no clear criterion is known which would determine the adequate class of abelian groups $E\stackrel{\pi}{\longrightarrow}L$ for a given loop $L$. Note that once such class is chosen, methods of non-associative Lie theory provide us with the classes of extensions for the corresponding Hopf and Sabinin algebras.

\subsection{The Ado theorem}\label{section:ado}

The Ado theorem is usually stated as a result about representations. While this may be appropriate in the context of Lie algebras, in the more general setting of Sabinin algebras it is more convenient to use a different wording which does not mention representations at all. 

Each Sabinin algebra appears as a Sabinin sublagebra of $\SU(A)$ for some non-associative algebra $A$; for instance, one can take $A$ to be the universal enveloping algebra $U(\el)$). Let us say that a variety of Sabinin algebras satisfies the Ado theorem if for each finite-dimensional Sabinin algebra $\el$ in that variety $A$ can be chosen to be finite-dimensional. By the universal property of $U(\el)$, this is the same as to say that the $U(\el)$ has an ideal of finite codimension which contains no non-trivial primitive elements.

We have mentioned before that, while Lie, Malcev and nilpotent Sabinin algebras satisfy the Ado theorem, this is not the case for all Sabinin algebras, Lie triple systems being a counterexample. The universal enveloping algebra of a central simple Lie triple system of finite dimension has no proper right ideals  apart from the augmentation ideal \cite{MP2, PI3}. It would be interesting to find criteria for a variety of loops to satisfy the Ado theorem.

The Ado theorem can be also considered as a statement about loops rather than their tangent algebras. Let us say that a local loop $Q$ is locally linear if can be locally embedded as a subloop of of the local loop of invertible elements in some finite-dimensional algebra $A$. This is equivalent to saying that the Sabinin algebra of $Q$ is isomorphic to a subalgebra of $\SU(A)$, which is the property described by the Ado theorem. Since the Ado theorem fails for Lie triple systems, local Bruck loops are not locally linear in general. 

We should point out that local linearity is a weaker concept than global linearity. This has nothing to do with non-associativity since the difference already exists for groups. For Moufang loops there also are examples which illustrate this phenomenon. For instance, the seven-dimensional projective plane is a Moufang loop which locally embeds into the invertible elements of the octonions. This local embedding, however, cannot be globalized \cite{Sh2}.

\section{Discrete loops and Sabinin algebras}
The definition of a of a nilpotent loop, and, generally, of the lower central series of a loop, was given by R.H.~Bruck in \cite{Bruck}. Since then, the theory of nilpotent groups has made significant advances and it became increasingly clear that Bruck's definition does not provide an adequate analogue of the associative lower central series. Namely, one should expect that
\begin{itemize}
\item the successive quotients of the lower central series for a finitely generated loop are finitely generated abelian groups; 
\item the graded abelian group associated with the lower central series carries an algebraic structure similar to that of a Lie ring;
\item there exists a close relation to the dimension series.
\end{itemize}
A lower central series for loops satisfying all the above properties has been defined and studied in \cite{M1, M2, MP1} and we review this construction in this section. A closely related definition was also given by F. Lemieux, C. Moore and D. Th\'erien in \cite{LMT}. Essentially, the argument of \cite{LMT} consists in defining the dimension series for the free loops and then pushing it to arbitrary loops using the universal property of the free loops. Recently, M.~Hartl and B.~Loiseau defined higher commutators in arbitrary semi-abelian categories which lead in the case of loops to the same definition as the one discussed here, see \cite{HL, HVdL}. 

\subsection{The commutator-associator series}\label{section:dev}
Let $L$ be a loop. The commutator-associator filtration on $L$ is  defined in terms of commutators, associators and associator deviations. The {commutator} of two elements $a,b$ of $L$ is
$$[a,b]=(ab)/(ba)$$
and the {associator} of $a,b$ and $c$ is defined by
$$(a,b,c)=((ab)c)/(a(bc)).$$
There is an infinite number of {associator deviations}. These
are functions $L^{l+3}\to L$ characterized by a non-negative number
$l$, called {\em level} of the deviation, and $l$ indices
$\alpha_{1},\ldots,\alpha_{l}$ with $0<\alpha_{i}\leq i+2$. The
deviations of level one are
$$(a,b,c,d)_{1}=(ab,c,d)/((a,c,d)(b,c,d)),$$
$$(a,b,c,d)_{2}=(a,bc,d)/((a,b,d)(a,c,d)),$$
$$(a,b,c,d)_{3}=(a,b,cd)/((a,b,c)(a,b,d)).$$
By definition, the deviation
$(a_{1},\ldots,a_{l+3})_{\alpha_{1},\ldots,\alpha_{l}}$ of level $l$
is equal to
$$A(a_{\alpha_{l}}a_{\alpha_{l}+1})/(A(a_{\alpha_{l}})A(a_{\alpha_{l}+1})),$$ where $A(x)$ stands for the
deviation $(a_{1},\ldots,
a_{\alpha_{l}-1},x,a_{\alpha_{l}+2},\ldots,
a_{l+3})_{\alpha_{1},\ldots,\alpha_{l-1}}$ of level $l-1$. The
associator is thought of as the associator deviation of level zero.

Now, set $\gamma_{1}L=L$ and for $n>1$ define $\gamma_{n}L$
to be the minimal normal subloop of $L$ containing
\begin{itemize}
\item{$[\gamma_{p}L,\gamma_{q}L]$ with $p+q\geq n$;}
\item{$(\gamma_{p}L,\gamma_{q}L,\gamma_{r}L)$
with $p+q+r\geq n$;}
\item{$(\gamma_{p_{1}}L,\ldots,\gamma_{p_{l+3}}L)_{
\alpha_{1},\ldots,\alpha_{l}}$ with $p_{1}+\ldots+p_{l+3}\geq n$.}
\end{itemize}
The subloop $\gamma_{n}L$ is called the {$n$th
commutator-associator subloop} of $L$. For a group $G$ the subgroup
$\gamma_{n}G$ is the $n$th term of the lower central series of $G$.

The commutator-associator subloops of a loop $L$ are normal in $L$. Moreover, they are fully invariant, that is, are preserved by all automorphisms of $L$. If $L$ is finitely generated, each quotient $\gamma_{i}L/\gamma_{i+1}L$ is a finitely generated abelian group. The crucial property of the commutator-associator filtration is that for an arbitrary loop $L$ the commutator, the associator and the associator deviations induce multilinear operations on the graded abelian group 
$$\mathcal{L}_\gamma L=\bigoplus \gamma_{i}L/\gamma_{i+1}L;$$ 
these operations respect the grading. It can be shown that  
the algebraic structure given by these
multilinear operations on $\mathcal{L}_\gamma L \otimes \Q$ is precisely that of a Sabinin algebra. This, however, is quite non-trivial and requires an understanding of the relation between the commuator-associator filtration and the dimension series.

\subsection{The dimension series} Let $L$ be a discrete loop and $\Q L$ its loop algebra over the rational numbers. Denote by $\Delta\subset \Q L$ the augmentation ideal and by $\Delta^n$ its $n$th power, that is, the submodule of $\Q L$ spanned over $\Q$ by all products of at least $n$ elements of $\Delta$ with any arrangement of the brackets. The loop $L$ can be thought of a subset of $\Q L$ and we define the $n$th  dimension subloop of $L$ as
$$D_nL = L \cap (1+\Delta^n).$$
It can be shown that $D_n L$ is, indeed, a series of fully invariant subloops of $L$. The successive quotients $D_i L/ D_{i+1} L$ are torsion-free abelian groups, and it turns out that the commutator, the associator and the associator deviations respect the filtration by the $D_n L$. The induced operations on the associated graded group can be identified as follows.

The loop algebra $\Q L$ is a bialgebra, with the coproduct defined as $\delta(g)=g\otimes g$ for $g\in L$. The associated graded algebra
$$\mathcal{D} L = \bigoplus_{i\geq 0} \Delta^i/\Delta^{i+1}$$ is then a cocommutative non-associative primitively generated Hopf algebra. If we set 
$$\mathcal{L} L =\bigoplus_{i>0} D_i L/ D_{i+1} L, $$
then $\mathcal{L} L\otimes\Q\subset \mathcal{D} L$ and
the primitive elements of $\mathcal{D} L$ are precisely the elements of $\mathcal{L} L \otimes \Q.$ Therefore, $\mathcal{L} L\otimes\Q$ is a Sabinin algebra. It can be shown that the Shestakov-Umirbaev operations $p_{m,n}$ on $\mathcal{D} L$ coincide with the operations induced by certain associator deviations on $\mathcal{L} L\otimes \Q$. (See \cite{MP1} for a more precise statement.)

These results can be applied to the study of the commutator-associator filtration. For each $n$ the subloop $D_n L$ contains the commutator-associator subloop $\gamma_n L$. Moreover, the analogue of the Jennings theorem \cite{Jennings} is true: 
$$D_nL = \sqrt{\gamma_n L},$$
that is, $D_n L$ consists of all those elements of $L$ whose $k$th power, for at least one $k$ and at least one arrangement of the parentheses, lies in $\gamma_n L$. In particular, we have that
$$\mathcal{L} L \otimes\Q= \mathcal{L}_\gamma L\otimes \Q.$$
 
\begin{remark}
The dimension subloops can be defined with the help of the group algebra with coefficients in any ring, and the result may depend on 
this ring. A particularly interesting case is that of the integer coefficients. For many years it was conjectured that the integer dimension subgroups of a group coincide with its lower central series. The eventual counterexample published by E.~Rips in 1972 was very involved, and while by now there are more counterexamples to this conjecture, the nature of the difference between the integer dimension series and the lower central series is still absolutely mysterious. It might be interesting to find a loop where the difference between the integer dimension series and the lower central series is due to the non-associative effects. 
\end{remark}

\subsection{The Magnus map and the dimension series of a free loop}

Let $\Z\{\{ X_1,\ldots, X_n\}\}$ be the non-associative ring of formal power series in $n$ non-commutative and non-associative variables. The invertible elements in this ring are the power series which start with $\pm 1$. They form a loop, which we denote by $\Z\{\{  X_1,\ldots, X_n\}\}^*$.
Let $F_{ \{ n \} }$ the free loop  on $n$ generators $x_1,\ldots, x_n$. 
The homomorphism
\begin{align*}
\mathcal{M}: F_{ \{ n \} } &\to \Z\{\{ X_1,\ldots, X_n \}\}^*,\\
x_i&\mapsto 1+X_i,
\end{align*}
is the non-associative version of the map~(\ref{eq:magnus}). As in the associative case, the $i$th term of the dimension series of $F_{ \{ n \} }$ consists of those elements which are sent by $\mathcal{M}$ to the power series of the form
$$1+\text{terms of degree at least\ } i.$$
This, however, is not enough to conclude that the loop $F_{ \{ n \} }$ is residually nilpotent since it is not known at the moment whether $\mathcal{M}$ is injective, even in the case of the free loop on one generator. 

A similar construction can be performed for the free commutative loop. It not known if the corresponding Magnus map is injective in this case, either.

The Magnus map identifies the completion of the Sabinin algebra $\mathcal{L} F_{ \{ n \} } \otimes\Q$ with the primitive elements in $\Z\{\{ X_1,\ldots, X_n\}\}$. In particular, this Sabinin algebra is the free Sabinin algebra on $n$ generators.

\section{Quantum loops}

An analogy with associative Hopf algebras suggests that there might exist deformations of the universal enveloping algebras for Sabinin algebras which are not necessarily coassociative or cocommutative. So far, no interesting examples of this kind have been found. Let us indicate, nevertheless, what one may be looking for.

In view of the fact that the Hopf algebra of distributions supported at the identity of a loop naturally satisfies the linearizations of the identities satisfied by the loop, one may assume (rightly or wrongly) that its possible deformations should satisfy the same linearized identities. Let us illustrate this logic with the variety of Moufang loops. 

The left and right Moufang identities can be written as
\begin{displaymath}
\gbeg{4}{9}
\got{2}{x} \raisebox{-0.07cm}{\got{1}{y}} \got{1}{\ }  \got{1}{z}    \gnl
\gcmu       \gcl{1}     \gcl{1}     \gnl
\gcl{1}     \gbr        \gcl{1}     \gnl
\gcl{1}     \gcl{1}     \gmu        \gnl
\gcl{1}     \gcl{1}     \gcn{1}{1}{2}{1}    \gnl
\gcl{1}     \gmu                    \gnl
\gcl{1}     \gcn{1}{1}{2}{1}        \gnl
\gmu                                \gnl
                                    \gnl
\got{1}{x(y(xz))}                   \gnl
\gend
=
\gbeg{4}{9}
\got{2}{x} \raisebox{-0.07cm}{\got{1}{y}} \got{1}{\ } \got{1}{z}    \gnl
\gcmu       \linea      \linea      \gnl
\linea      \cruce      \linea      \gnl
\gmu        \linea      \linea      \gnl
\gcn{2}{1}{2}{3}        \linea      \linea  \gnl
\gvac{1}    \gmu        \linea      \gnl
\gvac{1}    \gcn{2}{1}{2}{3}        \linea  \gnl
\gvac{2}    \gmu         \gnl
                                    \gnl
\gvac{2}     \got{1}{((xy)x)z}      \gnl
\gend
\quad\quad\mbox{and}\quad\quad
\gbeg{4}{9}
\got{1}{x} \raisebox{-0.07cm}{\got{1}{y}} \got{1}{\ } \got{2}{z}    \gnl
\linea \linea \gcmu \gnl
\linea \cruce \linea \gnl
\linea \linea \gmu \gnl
\linea \linea \gcn{1}{1}{2}{1}\gnl
\linea \gmu \gnl
\linea \gcn{1}{1}{2}{1}\gnl
\gmu \gnl
 \gnl
\got{1}{x(z(yz))} \gnl
\gend
=
\gbeg{4}{9}
\got{1}{x} \raisebox{-0.07cm}{\got{1}{y}} \got{1}{\ } \got{2}{z}    \gnl
\linea \linea \gcmu \gnl
\linea \cruce \linea \gnl
\gmu \linea \linea \gnl
\gcn{2}{1}{2}{3} \linea \linea \gnl
\gvac{1}\gmu \linea \gnl
\gcn{3}{1}{4}{5} \linea\gnl
\gvac{2}\gmu \gnl
\gnl
\gvac{2}\got{1}{((xz)y)z}\gnl
\gend
\end{displaymath}
where
$$\gbeg{2}{3}
    \got{1}{x} \raisebox{-0.07cm}{\got{1}{y}}   \gnl
    \gmu                    \gnl
    \gend = xy
\text{\ and\ }
\gbeg{2}{3}
    \got{2}{x}\gnl
    \gcmu                   \gnl
    \gend = \delta(x)$$
with  $\delta(x) = (x,x)$. The same diagrams represent the Moufang-Hopf identities satisfied in the corresponding algebra of distributions, with the difference that in this case $\delta(x) = \sum x_{(1)} \otimes x_{(2)}$. One may suspect that the universal enveloping algebra of the central simple  exceptional Malcev algebra of traceless octonions $\mathbb{M}(\alpha,\beta,\gamma)$ has non-cocommutative deformations that satisfy the Moufang-Hopf identities. Given that Hopf algebras are, in spirit, self-dual objects, it is natural to require that these deformations satisfy not coassociativity but, rather, the identities dual to Moufang-Hopf identities. These latter can be represented as 
$$
    \mbox{(left co-Moufang-Hopf)}
        \gbeg{4}{9}
        \gcmu \gnl
        \gcl{1} \gcn{1}{1}{1}{2}        \gnl
        \gcl{1} \gcmu                   \gnl
        \gcl{1} \gcl{1} \gcn{1}{1}{1}{2}\gnl
        \gcl{1} \gcl{1} \gcmu           \gnl
        \gcl{1} \gbr    \gcl{1}         \gnl
        \gmu    \gcl{1} \gcl{1}         \gnl
        \gend
        =
        \gbeg{4}{9}
        \gvac{2}    \gcmu       \gnl
        \gcn{3}{1}{5}{4}        \gcl{1}     \gnl
        \gvac{1}    \gcmu       \gcl{1}     \gnl
        \gcn{2}{1}{3}{2}        \gcl{1}     \gcl{1} \gnl
        \gcmu       \gcl{1}     \gcl{1}     \gnl
        \gcl{1}     \gbr        \gcl{1}     \gnl
        \gmu        \gcl{1}     \gcl{1}     \gnl
        \gend
        \quad\quad
        \mbox{(right co-Moufang-Hopf)}
        \gbeg{4}{9}
        \gcmu \gnl
        \gcl{1} \gcn{1}{1}{1}{2}        \gnl
        \gcl{1} \gcmu                   \gnl
        \gcl{1} \gcl{1} \gcn{1}{1}{1}{2}\gnl
        \gcl{1} \gcl{1} \gcmu           \gnl
        \gcl{1} \gbr    \gcl{1}         \gnl
        \gcl{1} \gcl{1} \gmu            \gnl
        \gend
        =
        \gbeg{4}{9}
        \gvac{2}    \gcmu       \gnl
        \gcn{3}{1}{5}{4}        \gcl{1}     \gnl
        \gvac{1}    \gcmu       \gcl{1}     \gnl
        \gcn{2}{1}{3}{2}        \gcl{1}     \gcl{1} \gnl
        \gcmu       \gcl{1}     \gcl{1}     \gnl
        \gcl{1}     \gbr        \gcl{1}     \gnl
        \gcl{1}     \gcl{1}     \gmu        \gnl
        \gend
$$
These identities can be written in the form of equations as 
$$
      \sum x_{(1)}x_{(2)(2)(1)} \, {\otimes}\,  x_{(2)(1)}\, {\otimes} \, x_{(2)(2)(2)} =
      \sum x_{(1)(1)(1)} x_{(1)(2)}\,  {\otimes}\, x_{(1)(1)(2)}\, {\otimes}\, x_{(2)}
$$
and
$$
      \sum x_{(1)}\, {\otimes}\, x_{(2)(2)(1)}\, {\otimes}\, x_{(2)(1)} x_{(2)(2)(2)} =
      \sum x_{(1)(1)(1)}\, {\otimes}\, x_{(1)(2)}\, {\otimes}\, x_{(1)(1)(2)}x_{(2)}.
$$
They arise in the study of the algebraic seven-dimensional sphere \cite{KM}.

Surprisingly, $U(\mathbb{M}(\alpha,\beta,\gamma))$ turns out to be rigid in the sense that any deformation which satisfies Moufang-Hopf and co-Moufang-Hopf identities is equivalent to a trivial deformation. Even universal enveloping algebras of finite-dimensional central simple Lie algebras exhibit certain rigidity if considered as Moufang-Hopf algebras since all their deformations are necessarily associative and coassociative; in other words, are quantized enveloping algebras in the usual sense \cite{MaP1,ShPI2}.

The rigidity characteristic of non-associativity is still poorly understood and needs an adecuate cohomological interpretation.

\begin{remark}
There are a number of deformations of non-associative algebras that carry the adjective ``quantum'', such as the quantum octonions of \cite{BeP} and of \cite{Bre-1}. They do not seem to fit in the framework of hypothetical quantum loops that we discuss here.
\end{remark}

\subsection*{Acknowledgments}
The authors were supported by the following grants:  FAPESP 2012/21938-4 and CONACyT grant 168093-F (J.M.), FAPESP 2012/22537-3 (J.M.P.-I.), CNPq 3305344/2009-9 and FAPESP 2010/50347-9 (I.P.Sh.), Spanish Government project MTM 2010-18370-C04-03 (all three). The first two authors would like to thank the Institute of Mathematics and Statistics of the University of S\~ao Paulo for hospitality.

\end{document}